\pdfoutput=1

\documentclass[
	english
]{jnsao}

\usepackage[utf8]{inputenc}

\usepackage{preamble}

\usepackage{enumitem}
\setenumerate{label=(\roman*)}

\numberwithin{equation}{section} 

\usepackage[capitalise,nameinlink]{cleveref}

\usepackage{mathtools}

\providecommand\given{\nonscript\;\delimsize|\nonscript\;\mathopen{}}
\DeclarePairedDelimiterX\set[1]\{\}{#1}
\DeclarePairedDelimiterX\seq[1](){#1}

\DeclarePairedDelimiterX\dual[2]{\langle}{\rangle}{#1,#2}

\DeclarePairedDelimiter\abs{\lvert}{\rvert}
\DeclarePairedDelimiter\norm{\lVert}{\rVert}

\DeclarePairedDelimiter\parens()
\DeclarePairedDelimiter\bracks[]

\newcommand\N{\mathbb{N}}
\newcommand\R{\mathbb{R}}
\newcommand\barR{\bar\R}

\renewcommand\d{\mathop{}\!\mathrm{d}}
\newcommand\dx{\d x}

\newcommand{\weaklystar}{\stackrel\star\rightharpoonup}

\newcommand{\dualspace}{^\star}
\newcommand{\bidualspace}{^{\star\star}}
\newcommand{\conjugate}{^\star}

\newcommand\BB{\mathcal{B}}
\newcommand\HH{\mathcal{H}}
\newcommand\KK{\mathcal{K}}
\newcommand\MM{\mathcal{M}}
\newcommand\ZZ{\mathcal{Z}}

\newcommand\Uad{U_{\mathrm{ad}}}

\DeclareMathAlphabet{\mathpzc}{OT1}{pzc}{m}{it}
\newcommand\oo{\mathpzc{o}}

\newcommand\sign{\operatorname{sign}}
\newcommand\dom{\operatorname{dom}}
\newcommand\supp{\operatorname{supp}}
\newcommand\interior{\operatorname{int}}

\let\subseteq\subset

\let\epsilon\varepsilon

\theoremstyle{definition}
\newtheorem{assumption}[theorem]{Assumption}
\crefname{assumption}{Assumption}{Assumptions}

\title{Second-order conditions for non-uniformly convex integrands: quadratic growth in \texorpdfstring{$\scriptstyle L^1$}{L¹}}
\shorttitle{Second-order conditions: quadratic growth in $L^1$}

\author{%
	Daniel Wachsmuth%
	\thanks{%
		Universität Würzburg,
		Institut für Mathematik,
		97074 Würzburg, Germany,
		\email{daniel.wachsmuth@mathematik.uni-wuerzburg.de},
		\url{https://www.mathematik.uni-wuerzburg.de/optimalcontrol}
	}
	\and
	Gerd Wachsmuth%
	\thanks{%
		Brandenburgische Technische Universität Cottbus--Senftenberg,
		Institute of Mathematics,
		03046 Cottbus,
		Germany,
		\email{gerd.wachsmuth@b-tu.de},
		\url{https://www.b-tu.de/fg-optimale-steuerung},
		\orcid{0000-0002-3098-1503}%
	}
}
\shortauthor{Wachsmuth, Wachsmuth}

\date{2022-03-09}

\acknowledgements{
This research was supported by the German Research Foundation (DFG) under grant numbers WA 3626/3-2 and WA 3636/4-2
within the priority program ``Non-smooth and Complementarity-based Distributed Parameter
Systems: Simulation and Hierarchical Optimization'' (SPP 1962).%
}

\manuscriptlicense{CC-BY-NC-SA 4.0}

\manuscriptsubmitted{2021-11-22}
\manuscriptaccepted{2022-05-09}
\manuscriptvolume{3}
\manuscriptnumber{8733}
\manuscriptyear{2022}
\manuscriptdoi{10.46298/jnsao-2022-8733}

\begin{document}

\maketitle

\begin{abstract}
We study no-gap second-order optimality conditions for a non-uniformly convex and non-smooth integral functional.
The integral functional is extended to the space of measures. The obtained second-order derivatives
contain integrals on lower-dimensional manifolds.
The proofs utilize the convex pre-conjugate, which is an integral functional on the space of continuous functions.
Applications to non-smooth optimal control problems are given.
\end{abstract}


\begin{msc}
	\mscLink{49J53}, 
	\mscLink{49J52}, 
	\mscLink{49K30} 
\end{msc}

\section{Introduction}
\label{sec:introduction}
We consider optimization problems of the form
\begin{equation}
	\label{eq:problem}
	\text{Minimize}
	\quad
	F(u) + G(u)
	\qquad\text{w.r.t.\ } u \in L^1(\lambda),
\end{equation}
where $F \colon \dom(G) \to \R$
is assumed to be smooth
and
$G(u) = \int_\Omega g(u) \d\lambda$
for some convex but possibly non-smooth
$g \colon \R \to \bar\R := \R \cup \set{\infty}$.
Here, $\Omega \subset \R^d$ is assumed to be non-empty, open, and bounded
and $\lambda$ is the Lebesgue measure.

We are interested in deriving no-gap optimality conditions of second order
in cases where $g$ is not uniformly convex.
Here, the meaning of no-gap second-order optimality conditions is: the difference between necessary and
sufficient conditions is as small as in the finite-dimensional case, i.e., positive semi-definiteness vs.\ positive definiteness of second derivatives.

In this article,
we follow the approach by
\cite{ChristofWachsmuth2017:1}
and
interpret the problem in the space of measures
$\MM(\Omega) = C_0(\Omega)\dualspace$.
Under certain assumptions on the problem data
and on the possible minimizer $\bar u \in L^1(\lambda)$,
we obtain that the second-order condition
\begin{equation*}
	F''(\bar u) \mu^2
	+
	G''(\bar u, -F'(\bar u); \mu)
	\ge 0
	\qquad\forall \mu \in \MM(\Omega)
\end{equation*}
is necessary for local optimality, while the second-order condition
\begin{equation*}
	F''(\bar u) \mu^2
	+
	G''(\bar u, -F'(\bar u); \mu)
	>0
	\qquad\forall \mu \in \MM(\Omega) \setminus\{0\}
\end{equation*}
is sufficient for local optimality (the latter with quadratic growth with respect to the $L^1(\lambda)$-norm).
Here, $F''(\bar u)$ is a weak-$\star$ continuous and quadratic form which allows for a second-order Taylor-like expansion of $F$
while $G''(\bar u, -F'(\bar u); \cdot)$
is the so-called second subderivative of $G$.
The focus of this article is two-fold: (1) prove these second-order conditions and (2) derive formulas for $G''(\bar u, -F'(\bar u); \cdot)$.

Let us put our work in perspective.
We restrict the discussion to contributions
to second-order conditions for optimal control
problems in which the control lives in a $d$-dimensional set
with $d \ge 2$, since the one-dimensional case is significantly simpler.

		First of all,
		the classical, smooth theory,
		see, e.g.,
		\cite{CasasTroeltzsch2012,Casas2015},
		is not applicable.
		On the one hand, $G$ might fail to be twice differentiable.
		On the other hand,
		even if $G$ is twice differentiable,
		the second derivative $G''(\bar u)$
		fails to be a Legendre form
		since $G$ is not assumed to be uniformly convex.
		However, such an assumption is crucial,
		see, e.g., \cite[Eq.~(2.5)]{CasasTroeltzsch2012}.

		The first results
		for sufficient second-order conditions
		for problems with bang-bang solution structure
		can be found in
		\cite{Casas2012:1}.
		Therein, the author proves
		quadratic growth
		w.r.t.\ the
		difference of the states
		by assuming some coercivity on an extended critical cone.
		No necessary conditions of second order are given.

		In
		\cite{CasasWachsmuthsBangBang},
		the authors
		use a structural assumption on the growth of the adjoint state $\bar\varphi$
		and this allows to formulate a
		second-order sufficient condition on measures
		living on the set $\set{\bar\varphi = 0}$.
		This implies quadratic growth in
		$L^1(\lambda)$.
		Again, no second-order necessary conditions are given.

		The first second-order necessary conditions
		for problems with bang-bang structure are given
		in \cite{ChristofWachsmuth2017:1}.
		Again, the structural assumption was used
		and the problem was analyzed in spaces of measures.
		The main observation of this paper is that
		(the boundary of) the set
		$\Uad = \set{u \in L^\infty(\lambda) \given -1 \le u \le 1}$
		possesses ``curvature'' in the space $\MM(\Omega)$,
		although this set is polyhedric in all Lebesgue spaces $L^p(\lambda)$,
		which, in some sense, means that it does not possess curvature.
		By employing a curvature term in the second-order conditions,
		no-gap conditions were given
		which are equivalent to a quadratic growth condition
		in $L^1(\lambda)$.

		A main application for second-order conditions
		is the derivation of discretization error estimates.
		In the bang-bang context, this was first studied in
		\cite{CasasWachsmuthsBangBang2}.

Second-order derivatives for non-smooth, convex integral functionals were studied by \cite{Do} ($L^p$ with $1<p<\infty$)
and \cite{Levy1993} (including $L^1$).
By extending the integral functional to the space of measures
(and by using weak-$\star$ convergence and not strong convergence as in \cite{Levy1993}),
we obtain second-derivatives that are smaller than those given in these works.
Hence, second-order conditions based on our result are stronger than those based on earlier results.
A characteristic feature
of our results is the appearance of integrals on $(d-1)$-dimensional manifolds in the expression of $G''$.
This was first observed in \cite{ChristofMeyer2019}, see their second-order expansion of the $L^1$-norm on $H^1_0(\Omega)$ in \cite[Corollary 4.10]{ChristofMeyer2019}.
The interpretation of these additional lower-dimensional integrals as curvature in measure space is due to \cite{ChristofWachsmuth2017:1}, where $g$ was chosen to be
the indicator function of the interval $[-1,1]$.

A characteristic motive in our work is the use of the pre-conjugate $J\colon C_0(\Omega) \to \bar\R$ of $G$, i.e., $J\conjugate=G$.
We derive second-order derivatives of $J$, and identify the second-order derivative of $G$ as its convex conjugate.
Interestingly, all assumptions that led to these results are given in terms of $J$.

The paper is structured as follows.
We start by adapting the theory of \cite{ChristofWachsmuth2017:1}
to the setting of \eqref{eq:problem} in \cref{sec:second-order_conditions}.
Afterwards, we study second-order derivatives
of convex integral functionals on $C_0(\Omega)$
in \cref{sec:sod_integral_functionals_cont}.
These results are utilized in \cref{sec:sod_integral_functionals_measures}
to obtain the expressions for $G''$.
The main results are \cref{thm:SSC_wo2} (No-gap second-order conditions) and \cref{thm_main_1} (Strictly twice epi-differentiability  of $G$).
Two applications are given in \cref{sec:app}.

\section{No-Gap Second-Order Conditions}
\label{sec:second-order_conditions}
In this section, we consider the minimization problem
\begin{equation*}
	\label{eq:prob}
	\tag{P}
	\text{Minimize}\quad
	\Phi(x) := F(x) + G(x)
	\qquad\text{w.r.t. } x \in X.
\end{equation*}
Here, $G \colon X \to \bar\R = (-\infty, \infty]$
and $F \colon \dom(G) \to \R$
are given.
We are interested in necessary and sufficient conditions of second order,
such that the gap between both conditions is as small as in finite dimensions.
We transfer the results of \cite{ChristofWachsmuth2017:1},
in which the case $G = \delta_C$ for a set $C \subset X$
was considered.
The generalization to the above problem with a more general $G$
is rather straightforward.

Throughout this section, we always consider the following situation.
\begin{assumption}[standing assumptions and notation]~
	\label{asm:standing_assumption}
	\begin{enumerate}
		\item
			$X$ is the (topological) dual of a separable Banach space $Y$,
		\item $\bar x $ is a fixed element of the set $\dom(G)$ (the minimizer/candidate for a minimizer),
		\item
			\label{asm:standing_assumption:3}
			There exist $F'(\bar x) \in Y$ and a bounded bilinear form $F''(\bar x) \colon X \times X \to \R$ with
			\begin{equation}
			\label{eq:hadamard_taylor_expansion}
				\lim_{k \to \infty}
				\frac{F(\bar x + t_k h_k) - F(\bar x) - t_k F'(\bar x) h_k - \frac12 t_k^2 F''(\bar x) h_k^2}{t_k^2}
				=
				0
			\end{equation}
			for all sequences $\seq{h_k} \subset X$, $\seq{t_k} \subset \R^+ := (0,\infty)$
			satisfying $t_k \searrow 0$,  $h_k \weaklystar h \in X$
			and
			$\bar x + t_k h_k \in \dom(G)$.
	\end{enumerate}
\end{assumption}
Note that we use the abbreviations $ F'(\bar x) h  := \dual{F'(\bar x)}{h}$
and $F''(\bar x) h^2 := F''(\bar x) ( h,h )$ for all $h \in X$ in \eqref{eq:hadamard_taylor_expansion},
and that \eqref{eq:hadamard_taylor_expansion} is automatically satisfied if $F$ admits a second-order Taylor expansion of the form
\begin{equation}
	\label{eq:strong_taylor_expansion}
	F(\bar x + h)
	-
	F(\bar x)
	-
	F'(\bar x) h
	-
	\frac12 F''(\bar x) h^2
	=
	\oo(\norm{h}^2_X)
	\quad\text{ as } \norm{h}_X \to 0.
\end{equation}

As a second derivative for the functional $G$,
we use the so-called
weak-$\star$ second subderivative.
\begin{definition}[weak-\texorpdfstring{$\star$}{*} second subderivative]
	\label{def:weak_star_subderivative}
	Let $x \in \dom(G)$ and $w \in Y$ be given.
	Then the weak-$\star$ second subderivative
	$G''(x, w; \cdot) \colon X \to [-\infty, \infty]$
	of $G$ at $x$ for $w$ is defined by
	\begin{equation*}
		G''(x, w; h)
		:=
		\inf
		\set*{
			\liminf_{k \to \infty} \frac{G(x + t_k h_k) - G(x) - t_k \dual{w}{h_k}}{t_k^2/2}
			\given
			t_k \searrow 0,
			h_k \weaklystar h
		}
		.
	\end{equation*}
\end{definition}

	We discuss some properties of $G''(x,w;\cdot)$ for a convex function $G$
	and $x \in \dom(G)$.
	In the following, the parameter $w$ will often be taken from the convex subdifferential $\partial G(x)$.
	We recall that
	$\partial G(x)$ is a subset of the dual space $X\dualspace$.
	In the case that $X$ is not reflexive, $X\dualspace$ is bigger than $Y$
	(note that we identify $Y$ with a subspace of $X\dualspace = Y\bidualspace$ in the canonical way).
	Thus, the existence of $w\in Y \cap \partial G(x)$ is an additional regularity assumption, since
the existence of subgradients of $G$ in the smaller space $Y$
is not guaranteed in general,
as the next example demonstrates.

\begin{remark}
	\label{rem:no_predual_subgradients}
	We choose $Y = c_0$ (zero sequences equipped with supremum norm).
	Thus, the dual spaces are (isometric to)
	$X = Y\dualspace = \ell^1$ and $X\dualspace = \ell^\infty$.
	By using
	\begin{equation*}
		C := \set*{ x \in \ell^1 \given \forall n \in \N : \abs{x_n} \le n^{-2} }
	\end{equation*}
	we define $G \colon \ell^1 \to \barR$ via
	\begin{equation*}
		G(x) = \sum_{n=1}^\infty x_n \in \R \qquad\forall x \in C
	\end{equation*}
	and $G(x) = \infty$ for all $x \in \ell^1 \setminus C$.

	The set $C$ is bounded due to $\sum_{n = 1}^\infty 1/n^2 < \infty$
	and weak-$\star$ closed since it is the intersection of the weak-$\star$ closed ``stripes''
	\begin{equation*}
		\set*{x \in \ell^1 \given \abs{x_n} \le n^{-2} } \qquad\forall n \in \N.
	\end{equation*}
	Thus, $C$ is weak-$\star$ compact.
	The function $G$ is convex.
	In order to check that $G$ is weak-$\star$ continuous on $C$,
	let $x_0 \in C$ be given and consider a net $\seq{x_i}_{i\in I} \subset C$ with $x_i \to x_0$.
	For an arbitrary $\varepsilon > 0$, there is $N \in \N$ with
	$\sum_{n = N+1}^\infty n^{-2} < \varepsilon$.
	Next, there is $i \in I$ with
	\begin{equation*}
		\abs*{ \sum_{n = 1}^N (x_{j,n} - x_{0,n}) } < \varepsilon
		\qquad\forall j \ge i
	\end{equation*}
	since $y \mapsto \sum_{n = 1}^N y_n$ is weak-$\star$ continuous.
	Thus,
	\begin{equation*}
		\abs{G(x_j) - G(x_0)}
		\le
		\abs*{ \sum_{n = 1}^N (x_{j,n} - x_{0,n}) }
		+
		\sum_{n = N+1}^\infty \abs{x_{j,n}}
		+
		\sum_{n = N+1}^\infty \abs{x_{0,n}}
		<
		3 \varepsilon
		\qquad\forall j \ge i.
	\end{equation*}
	Since $\varepsilon > 0$ was arbitrary, this shows weak-$\star$ continuity on $C$.
	Hence, $G$ is weak-$\star$ lower semicontinuous on $\ell^1$.
	Finally, it is easy to check that $\partial G(0) = \set{1}$,
	but $1 \in \ell^\infty \setminus c_0$.

	A sufficient condition for the existence of subgradients in the predual space at $x$
would be the weak-$\star$ continuity of $G$ at $x$. In infinite dimensions, this is a very restrictive condition, since
it implies that the convex function $G$ is bounded (and hence constant) on a subspace with finite codimension.
\end{remark}

We continue with some basic properties of $G''(x, w; \cdot )$.

\begin{lemma}
	\label{lem:Gpp_convex}
	We assume that $G$ is convex and $x \in \dom(G)$.
	For $w \in Y \cap \partial G(x)$ we have
	\begin{equation*}
		\forall h \in X : G''(x, w; h ) \ge 0,
	\end{equation*}
	whereas in case $w \in Y \setminus \partial G(x)$ we have
	\begin{equation*}
		\exists h \in X \setminus \set{0} : G''(x, w; h ) = -\infty.
	\end{equation*}
\end{lemma}
\begin{proof}
	In case $w \in Y \cap \partial G(x)$, the definition of $\partial G(x)$ implies
	$G(x + t_k h_k) - G(x) - \dual{w}{t_k h_k} \ge 0$
	for arbitrary $t_k > 0$ and $h_k \in X$.
	This directly yields
	$G''(x, w; h ) \ge 0$ for all $h \in X$.

	In case $w \in Y \setminus \partial G(x)$, there exists
	$h \in X \setminus \set{0}$ such that
	$G(x+h)-G(x) = \dual{w}{h} - \tau$
	for some $\tau > 0$.
	We choose $\seq{t_k}_{k \in \N} \subset (0,1)$ with $t_k \searrow 0$.
	By convexity, we get
	\begin{equation*}
		\frac{G(x + t_k h) - G(x) - t_k \dual{w}{h}}{t_k^2/2}
		\le
		\frac{
			G(x + h) - G(x) - \dual{w}{h}
		}{t_k/2}
		\le
		-\frac2{t_k}\tau
		\to
		-\infty
		.
	\end{equation*}
	Thus, $G''(x, w; h) = -\infty$.
\end{proof}

	We recall that the convexity of $G$ implies
	$ G'(x; h) \ge \dual{w}{h} $
	for all $x \in \dom(G)$, $w \in Y \cap \partial G(x)$, and $h \in X$.
	The next lemma can be used to
	reduce the second-order conditions below
	to the so-called critical cone
	on which this inequality is satisfied with equality.
\begin{lemma}
	\label{lem:Gpp_directional_derivative}
	We assume that $G$ is convex and $x \in \dom(G)$
	is chosen such that the directional derivative
	$G'(x; \cdot) \colon X \to [-\infty, \infty]$
	is sequentially weak-$\star$ lower semicontinuous.
	Further let $w \in Y$ be given.
	Then
	\begin{equation*}
		G'(x; h) > \dual{w}{h}
		\;\Rightarrow\;
		G''(x, w; h) = + \infty
		\qquad\forall h \in X.
	\end{equation*}
\end{lemma}
\begin{proof}
	Let $h \in X$ with $G'(x; h) > \dual{w}{h}$ be arbitrary.
	For all sequences
	$h_k \weaklystar h$ and $t_k \searrow 0$
	we have
	$G(x + t_k h_k) - G(x) \ge t_k G'(x; h_k)$
	by convexity.
	Thus,
	\begin{align*}
		\liminf_{k \to \infty} \frac{G(x + t_k h_k) - G(x) - t_k \dual{w}{h_k}}{t_k}
		&\ge
		\liminf_{k \to \infty} \parens[\big]{G'(x; h_k) - \dual{w}{h_k}}
		\ge
		G'(x; h) - \dual{w}{h}
		>
		0
		.
	\end{align*}
	Since an additional factor $t_k^{-1}$ appears in the definition of $G''(x, w; h)$,
	this implies $G''(x,w; h) = +\infty$.
\end{proof}

The assumption on $G'$ is satisfied in the following situation:
Let $G$ be convex and lower semicontinuous.
	In the case that $X$ is reflexive,
	the weak topology and the weak-$\star$ topology coincide.
	Moreover, if $x \in \interior \dom(G)$
	we know that $G$ is locally Lipschitz at $x$.
	Consequently, $G'(x; \cdot)$ is convex and Lipschitz continuous,
	thus sequentially weak-$\star$ lower semicontinuous.
	In the non-reflexive situation,
	this argument is no longer applicable, and
	the statement of  \cref{lem:Gpp_directional_derivative} may fail
	for (strongly) continuous $G$
	as the following example shows.
\begin{remark}
	\label{rem:bad_non_reflexive}
Let $Y = C_0(\Omega)$, $X = Y\dualspace = \MM(\Omega)$
and define
$G(x) := \norm{x}_{\MM(\Omega)}$.
This functional is obviously convex, globally Lipschitz continuous,
and weak-$\star$ lower semicontinuous.
Let $K \subset \Omega$ be compact with nonempty interior.
We consider the measure
$x \in \MM(\Omega)$, $x(A) := \lambda( K \cap A)$,
and some $w \in C_0(\Omega)$ with $\norm{w}_{C_0(\Omega)} \le 1$ and $w = 1$ in $K$.
Thus,
$\dual{w}{x}  =\norm{x}_{\MM(\Omega)}$ and
$w \in Y \cap \partial G(x)$.
Next, we choose the Dirac measure $h = -\delta_\omega$ for some fixed $\omega \in \interior K$.
Clearly, $1 = G'(x; h) > \dual{w}{h} = -1$.
Finally, we choose $s_k \searrow 0$ and
define $t_k > 0$ and the measures $h_k$ via
\begin{equation*}
	t_k := \lambda( B_{s_k}(\omega)),
	\qquad
	h_k(A) :=
	-\frac{\lambda( B_{s_k}(\omega) \cap A )}{t_k}
	.
\end{equation*}
Note that $t_k \searrow 0$ and $h_k \weaklystar h$.
Then, we have (for $k$ large enough such that $B_{s_k}(\omega) \subset K$)
\begin{equation*}
	G(x + t_k h_k) - G(x) - t_k \dual{w}{h_k}
	=
	\int_K \abs*{1 - \chi_{B_{s_k}(\omega)}} \d\lambda
	-
	\int_K 1 \d\lambda
	+
	t_k
	=
	0.
\end{equation*}
Thus,
\begin{equation*}
	\frac{G(x + t_k h_k) - G(x) - t_k \dual{w}{h_k}}{t_k^2/2} \to 0
\end{equation*}
and this shows
$G''(x,w;h)=0$.
Note that $G'(x; \cdot)$ cannot be weak-$\star$ lower semicontinuous,
since this would contradict \cref{lem:Gpp_directional_derivative}.
\end{remark}
	In the next definition, we ensure the existence of recovery sequences.
\begin{definition}[second-order epi-differentiability]
	\label{def:second_order_derivative}
	Let $x \in  \dom(G)$ and $w \in Y$ be given.
	The functional $G$ is said to be weakly-$\star$ twice epi-differentiable
	(respectively, strictly twice epi-differentiable, respectively, strongly
	twice epi-differentiable)
	at $x$ for $w$ in a direction $h \in X$,
	if for all $\seq{t_k}\subset \R^+$ with $t_k \searrow 0$ there exists a sequence $\seq{h_k}$
	satisfying $h_k \weaklystar h$
	(respectively, $h_k \weaklystar h$ and $\norm{h_k}_X \to \norm{h}_X$, respectively, $h_k \to h$) and
	\begin{equation}
		\label{eq:recovery_sequence_def}
		G''(x, w; h )
		=
		\lim_{k \to \infty} \frac{G(x + t_k h_k) - G(x) - t_k \dual{w}{h_k}}{t_k^2/2}
		.
	\end{equation}
	The functional $G$ is called weakly-$\star$/strictly/strongly twice epi-differentiable at $x$ for $w$ if it is
	weakly-$\star$/strictly/strongly twice epi-differentiable at $x$ for $w$ in all directions $h \in X$.
\end{definition}

We provide second-order optimality conditions for \eqref{eq:prob}.
We start with the second-order necessary condition (SNC).

\begin{theorem}[SNC involving the second subderivative]
	\label{thm:SNC}
	Suppose that $\bar x$ is a local minimizer of \eqref{eq:prob} such that
	\begin{equation}
		\label{eq:second_order_growth}
		\Phi(x)
		\ge
		\Phi(\bar x) + \frac{c}{2} \norm{x - \bar x}^2_X
		\qquad\forall x \in B_\varepsilon^X(\bar x)
	\end{equation}
	holds for some $c \geq 0$ and some $\varepsilon > 0$.
	Here, $B_\varepsilon^X(\bar x)$ denotes the closed ball in $X$ centered at $\bar x$ with radius $\varepsilon$.
	Assume further that one of the following conditions is satisfied.
	\begin{enumerate}
		\item The map $h \mapsto F''(\bar x) h^2$ is sequentially weak-$\star$ upper semicontinuous. \label{assumption-usc}
		\item The functional $G$ is strongly twice epi-differentiable at $\bar x$ for $-F'(\bar x)$. \label{assumption-mrc}
	\end{enumerate}
	Then
	\begin{equation}
		\label{eq:weakstarSNC}
		F''(\bar x) h^2
		+
		G''(\bar x, -F'(\bar x); h)
		\ge
		c \norm{h}_X^2
		\qquad
		\forall h \in X.
	\end{equation}
\end{theorem}
\begin{proof}
	We first consider the case \ref{assumption-usc}:
	Let $h \in X$ be given.
	In case $G''(\bar x, -F'(\bar x); h) = \infty$,
	\eqref{eq:weakstarSNC} holds automatically.
	Otherwise,
	it follows
	from the definition of $G''(\bar x, -F'(\bar x); \cdot)$,
	that for every $M > G''(\bar x; -F'(\bar x); h)$ we can find sequences
	$\seq{h_k} \subset X$ and $\seq{t_k} \subset \R^+$ such that
	$t_k \searrow 0$,
	$h_k \weaklystar h$,
	$x_k := \bar x + t_k h_k$
	and
	\begin{equation}
		\label{eq:choice_of_rk}
		\lim_{k \to \infty} \frac{G(\bar x + t_k h_k) - G(\bar x) + t_k \dual{F'(\bar x)}{h_k}}{t_k^2/2} \le M.
	\end{equation}
	Since $x_k \to \bar x$ strongly in $X$, \eqref{eq:second_order_growth} entails
	$\Phi(x_k) \ge \Phi(\bar x) + \frac{c}{2} \norm{x_k - \bar x}^2_X $ for  large enough $k$.
	Adding \eqref{eq:hadamard_taylor_expansion} and \eqref{eq:choice_of_rk}
	yields
	\begin{align*}
		M
		&\ge
		\lim_{k \to \infty} \frac{\Phi(\bar x + t_k h_k) - \Phi(\bar x) - \frac12 t_k^2 F''(\bar x) h_k^2}{t_k^2/2}
		\ge
		\limsup_{k \to \infty} \parens*{ c \norm{h_k}^2 - F''(\bar x) h_k^2 }
		\\
		&\ge
		\liminf_{k \to \infty} \parens*{ c \norm{h_k}^2} + \limsup_{k \to \infty} \parens*{ - F''(\bar x) h_k^2 }
		=
		\liminf_{k \to \infty} \parens*{ c \norm{h_k}^2} - \liminf_{k \to \infty} \parens*{ F''(\bar x) h_k^2 }
		\\
		&\ge
		c \norm{h}^2 - F''(\bar x) h^2
		.
	\end{align*}
	Here, we used that the function $h \mapsto F''(\bar x) h^2$ is sequentially weak-$\star$ upper semicontinuous.
	Since $M > G''(\bar x, -F'(\bar x); h)$ was arbitrary,
	we obtain \eqref{eq:weakstarSNC}.

	It remains to prove \eqref{eq:weakstarSNC} under assumption \ref{assumption-mrc}.
	For every $h \in X$
	we can find sequences
	$\seq{h_k} \subset X$, $\seq{t_k} \subset \R^+$ such that
	$t_k \searrow 0$,
	$h_k \to h$,
	$x_k := \bar x + t_k h_k$
	and
	\begin{equation*}
		\lim_{k \to \infty} \frac{G(\bar x + t_k h_k) - G(\bar x) + t_k \dual{F'(\bar x)}{h_k}}{t_k^2/2}
		=
		G''(\bar x, -F'(\bar x); h)
		.
	\end{equation*}
	Now,
	the second-order condition \eqref{eq:weakstarSNC} follows
	analogously to case \ref{assumption-usc}.
	Note that here $h_k\to h$ implies $F''(\bar x) h_k^2\to F''(\bar x) h^2$.
\end{proof}

We continue with the second-order sufficient condition (SSC).

\begin{theorem}[SSC involving the second subderivative]
	\label{thm:SSC_wo}
	Assume that the map $h \mapsto F''(\bar x) h^2$
	is sequentially weak-$\star$ lower semicontinuous
	and that
	\begin{equation}
		\label{eq:SSC_wo}
		F''(\bar x) h^2
		+
		G''(\bar x, -F'(\bar x); h)
		>
		0
		\qquad
		\forall h \in X \setminus \set{0}.
	\end{equation}
	Suppose further that
	\begin{equation*}
		\label{eq:NDS}
		\tag{\textup{NDC}}
		\begin{aligned}[t]
			&\text{for all $\seq{h_k} \subset X$, $\seq{t_k} \subset \R^+$
			with $h_k \weaklystar 0$, $t_k \searrow 0$
			and $\norm{h_k}_X = 1$, it is true that }\\
			&\qquad
			\liminf_{k \to \infty} \parens[\bigg]{
				\frac1{t_k^2} \parens[\big]{G(\bar x + t_k h_k) - G(\bar x)}
				+
				\dual{F'(\bar x)}{h_k / t_k}
				+
				\frac12 F''(\bar x) h_k^2
			}
			> 0
			.
		\end{aligned}
	\end{equation*}
	Then $\bar x$ satisfies the growth condition \eqref{eq:second_order_growth} with some constants $c > 0$ and $\varepsilon > 0$.
\end{theorem}
	In case $G$ is convex,
	a first-order condition is hidden in \eqref{eq:SSC_wo},
	see \cref{lem:Gpp_convex} and \cref{cor:2nd_order_traditional} below.
The acronym \eqref{eq:NDS}
stands for ``non-degeneracy condition'', compare \cite[Theorem 4.4]{ChristofWachsmuth2017:1}.
We will give a sufficient condition for \eqref{eq:NDS} below in \cref{lem:suf_nds}.

\begin{proof}
	We do not argue by contradiction.
	For a sequence $\seq{\varepsilon_k} \subset (0,\infty)$ with $\varepsilon_k \searrow 0$,
	we define
	\begin{equation*}
		c_k :=
		\inf\set[\bigg]{
			\frac{\Phi(x) - \Phi(\bar x)}{\frac12 \norm{x - \bar x}_X^2}
			\given
			x \in B_{\varepsilon_k}^X(\bar x) \setminus \set{\bar x}
		}
		\subset [-\infty, \infty]
		.
	\end{equation*}
	We have to show that $c_k > 0$ for some $k \in \N$.
	Note that $c_k$ is increasing
	and it is sufficient
	to prove that
	$\hat c := \lim_{k \to \infty} c_k > 0$.
	By definition of $c_k$,
	we can find a sequence $\seq{x_k}$ with $x_k \to \bar x$ in $X$
	and
	\begin{equation*}
		\hat c
		=
		\lim_{k \to \infty}
		\frac{\Phi(x_k) - \Phi(\bar x)}{\frac12 \norm{x_k - \bar x}_X^2}
		.
	\end{equation*}
	Define $t_k := \norm{x_k - \bar x}_X$
	and
	$h_k := (x_k - \bar x) / t_k$.
	Then $\norm{h_k}_X = 1$ for all $k$ and we may extract a subsequence (not relabeled)
	such that $h_k \weaklystar h$.
	Now, we utilize \eqref{eq:hadamard_taylor_expansion}
	and have
	\begin{align*}
		\hat c
		&=
		\lim_{k \to \infty}
		\frac{F(\bar x + t_k h_k) - F(\bar x) + G(\bar x + t_k h_k) - G(\bar x)}{t_k^2/2}
		\\
		&=
		\lim_{k \to \infty}
		\parens[\bigg]{
			\frac{G(\bar x + t_k h_k) - G(\bar x) + t_k F'(\bar x) h_k}{t_k^2/2}
			+
			F''(\bar x) h_k^2
		}
		.
	\end{align*}
	In case that $h = 0$, the condition
	\eqref{eq:NDS}
	implies that the right-hand side is positive and we are done.
	Otherwise, $h \ne 0$
	and we can use the definition of $G''$ together with the sequential weak-$\star$ lower semicontinuity of $F''(\bar x)$.
	This leads to
	\begin{equation*}
		\hat c
		\ge
		G''(\bar x, -F'(\bar x); h) + F''(\bar x) h^2.
	\end{equation*}
	The right-hand side is
	positive
	by \eqref{eq:SSC_wo}
	and this finishes the proof.
\end{proof}

By combining the previous two theorems,
we arrive at our main theorem on no-gap second-order conditions.
\begin{theorem}[no-gap second-order optimality condition]
	\label{thm:no-gap-SOC}
	Assume that the map $h \mapsto F''(\bar x) h^2$
	is sequentially weak-$\star$ lower semicontinuous,
	that  \eqref{eq:NDS} holds,
	and that one of the conditions \ref{assumption-usc} and \ref{assumption-mrc} in \cref{thm:SNC}
	is satisfied.
	Then, the condition
	\begin{equation*}
		F''(\bar x) h^2
		+
		G''(\bar x, -F'(\bar x); h)
		>
		0
		\qquad
		\forall h \in X \setminus \set{0}
	\end{equation*}
	is equivalent to the quadratic growth condition \eqref{eq:second_order_growth} with constants $c>0$ and $\varepsilon > 0$.
\end{theorem}
In the case that $G$ is convex,
we
can combine \cref{thm:no-gap-SOC} with \cref{lem:Gpp_convex,lem:Gpp_directional_derivative}.
Thus,
we can recast the above second-order conditions
in a familiar form including
the first-order condition and a critical cone.
\begin{corollary}
	\label{cor:2nd_order_traditional}
	In addition to the assumptions of \cref{thm:no-gap-SOC},
	we assume that $G$ is convex and that $G'(\bar x; \cdot)$
	is sequentially weak-$\star$ lower semicontinuous.
	The quadratic growth condition \eqref{eq:second_order_growth}
	with constants $c > 0$ and $\varepsilon > 0$ is satisfied
	if and only if
	\begin{subequations}
		\label{eq:2nd_order_traditional}
		\begin{align}
			\label{eq:2nd_order_traditional_first_order}
			F'(\bar x) + \partial G(\bar x) &\ni 0 \\
			\label{eq:2nd_order_traditional_second_order}
			F''(\bar x) h^2 + G''(\bar x, -F'(\bar x); h) &> 0 \qquad \forall h \in \KK \setminus \set{0},
		\end{align}
	\end{subequations}
	where
	the critical cone $\KK$ is defined via
	\begin{equation*}
		\KK := \set{ h \in X \given F'(\bar x) h + G'(\bar x; h) = 0 }.
	\end{equation*}
	Moreover, if $\bar x$ is a local minimizer of \eqref{eq:prob},
	then \eqref{eq:2nd_order_traditional} holds with
	``$\;\!\ge\!$'' instead of ``$\;\!>\!$''.
\end{corollary}
\begin{proof}
	Let \eqref{eq:second_order_growth} be satisfied
	with $c \ge 0$ and $\varepsilon > 0$.
	From \cref{thm:SNC}, we get \eqref{eq:weakstarSNC}.
	Clearly, \eqref{eq:2nd_order_traditional_second_order}
	(with ``$\ge$'' instead of ``$>$'' in case $c = 0$)
	follows.
	If \eqref{eq:2nd_order_traditional_first_order} is violated,
	\cref{lem:Gpp_convex} yields the existence of
	$h \in X \setminus \set{0}$ with $G''(x, w; h) = -\infty$
	and this contradicts \eqref{eq:weakstarSNC}.
	This shows the ``only if'' part of the first assertion
	and the second assertion.

	Let \eqref{eq:2nd_order_traditional} be satisfied.
	For $h \in X \setminus \set{0}$,
	we have
	$G'(\bar x; h) \ge -\dual{F'(\bar x)}{h}$
	due to \eqref{eq:2nd_order_traditional_first_order}.
	In case $h \in \KK \setminus \set{0}$, \eqref{eq:SSC_wo}
	follows from \eqref{eq:2nd_order_traditional_second_order}.
	Otherwise, we have
	$G'(\bar x; h) > -\dual{F'(\bar x)}{h}$
	and
	\eqref{eq:SSC_wo} is implied
	by \cref{lem:Gpp_directional_derivative}.
\end{proof}
If \eqref{eq:2nd_order_traditional_first_order} holds,
we have
$\KK = \set{ h \in X \given F'(\bar x) h + G'(\bar x; h) \le 0 }$.
Since $G'(\bar x; \cdot)$ is convex,
this results in the convexity of $\KK$.
In case $G'(\bar x; \cdot)$ is
(sequentially) weak-$\star$ lower semicontinuous,
$\KK$ is also (sequentially) weak-$\star$ closed.

Finally, we give a sufficient condition for \eqref{eq:NDS}.
\begin{lemma}
	\label{lem:suf_nds}
	Suppose that $h \mapsto F''(\bar x) h^2$ is sequentially weak-$\star$ lower semicontinuous
	and that there exist $c, \varepsilon > 0$ such that
	\begin{equation}
		\label{eq:suf_nds}
		G(x) - G(\bar x) + F'(\bar x) (x - \bar x)
		\ge
		\frac c2 \norm{ x - \bar x}_X^2
		\qquad
		\forall x \in B^X_\varepsilon(\bar x).
	\end{equation}
	Then, \eqref{eq:NDS} is satisfied.
\end{lemma}
\begin{proof}
	Let sequences as in \eqref{eq:NDS} be given.
	Then,
	\begin{align*}
		&\liminf_{k \to \infty} \parens[\bigg]{
			\frac1{t_k^2} \parens[\big]{G(\bar x + t_k h_k) - G(\bar x)}
			+
			\dual{F'(\bar x)}{h_k / t_k}
			+
			\frac12 F''(\bar x) h_k^2
		}
		\\&\quad
		\ge
		\liminf_{k \to \infty} \parens[\bigg]{
			\frac c2
			+
			\frac12 F''(\bar x) h_k^2
		}
		\ge
		\frac c2 + F''(\bar x) 0^2
		=
		\frac c2 > 0.
		\qedhere
	\end{align*}
\end{proof}

\section{Second-order derivatives of integral functionals over continuous functions}
\label{sec:sod_integral_functionals_cont}
Our goal is to apply the theory of \cref{sec:second-order_conditions}
in the setting $X = \MM(\Omega)$ and $Y = C_0(\Omega)$.
We will see in \cref{sec:sod_integral_functionals_measures}
that second subderivatives of integral functionals
on $\MM(\Omega)$
can be obtained by first studying the pre-conjugate functionals.
Therefore,
this section is devoted
to
integral functionals
\begin{equation}
	\label{eq:integral_functionals_cont}
	J(w)
	:=
	\int_\Omega j(w(\omega)) \d\lambda(\omega)
	\qquad
	\forall w \in C_0(\Omega),
\end{equation}
where $j \colon \R \to \R$ is a convex function and
$\Omega\subset\R^d$ is a bounded open set.
We equip $\Omega$ with the Borel $\sigma$-algebra $\BB(\Omega)$
and the $d$-dimensional Lebesgue measure $\lambda := \lambda^d$ (restricted to $\Omega$).
Note that $J(w) \in \R$ for all $w \in C_0(\Omega)$.
Later, we will choose $J$ such that its convex conjugate is equal to $G$.

The main goal of this section is to provide results
for the functional $J$ that allow to compute
the weak-$\star$ second subderivative (in the sense of \cref{def:weak_star_subderivative})
of the convex conjugate $G=J\conjugate$, which is a mapping from $\MM(\Omega) = C_0(\Omega)\dualspace$ to $\bar\R$.
In particular,
we are going to compute second-order derivatives of $J$
in two different ways.
First, we compute second subderivatives similar to \cref{def:weak_star_subderivative}
(but w.r.t.\ the strong topology)
in \cref{subsec:SOSD}
and, second,
we investigate difference quotients of subdifferentials
in \cref{subsec:dq_of_subdiff}.
We will see that (under appropriate assumptions)
both approaches
lead to the same object.
These findings will be crucial in \cref{sec:sod_integral_functionals_measures}
to prove that $G \colon \MM(\Omega) \to \bar \R$
is strictly twice epi-differentiable in the sense of \cref{def:second_order_derivative}.

The analysis in this section relies heavily on the following characterization of the
convex conjugate of integral functionals on continuous
functions for finite Carathéodory integrands.

\begin{theorem}[\texorpdfstring{\cite[Corollaries~4A, 4B]{Rockafellar1971}}{Rockafellar, 1971, Corollary~4A, 4B}]
	\label{thm:rockit}
	Let $\Omega \subset \R^d$ be
	open
	and
	let $\nu \colon \BB(\Omega) \to [0,\infty]$
	be a
	$\sigma$-finite,
	regular Borel measure.
	We suppose that $\iota \colon \Omega \times \R \to \R$
	is a convex Carathéodory function, i.e.,
	$\iota(\omega, \cdot)$ is convex and continuous for every $\omega \in \Omega$
	and
	$\iota(\cdot, x)$ is measurable for all $x \in \R$.
	Finally, we assume that $\iota(\cdot, x) \in L^1(\nu)$ for all $x \in \R$.
	Then, the integral functional $I \colon C_0(\Omega) \to \R$
	\begin{equation*}
		I(u)
		:=
		\int_\Omega \iota(\omega, u(\omega)) \d\nu(\omega)
		\qquad
		\forall u \in C_0(\Omega)
	\end{equation*}
	is well-defined, convex and continuous.
	Its convex conjugate $I\conjugate \colon \MM(\Omega) \to \bar\R$ is given by
	\begin{equation*}
		I\conjugate(\mu)
		=
		\int_\Omega \iota\conjugate\parens*{\omega , \frac{\d\mu}{\d\nu}(\omega)} \d\nu(\omega)
		\qquad
		\forall \mu \in \MM(\Omega), \mu \ll \nu
	\end{equation*}
	and
	$I\conjugate(\mu) = +\infty$ if $\mu$ is not absolutely continuous w.r.t.\ $\nu$,
	i.e., if $\mu \ll \nu$ does not hold.
	Here, $\iota\conjugate(\omega,\cdot)$ is the convex conjugate of $\iota(\omega,\cdot)$
	for every $\omega \in \Omega$.

	In addition, $\mu\in \partial I(u)$ if and only if $\mu \ll \nu$
	and
	\[
	 \frac{\d\mu}{\d\nu}(\omega) \in \partial \iota(\omega,u(\omega))
	 \quad\text{for $\nu$-almost all $\omega \in \Omega$.}
	\]
\end{theorem}
Note that the Fenchel--Young inequality implies that
$\iota\conjugate(\omega, \cdot) \ge -\iota(\omega, 0)$.
The right-hand side is integrable and thus, the integral
in the definition of $I\conjugate$ is well defined
(in the sense of quasi-integrability).

\subsection{Strong second subderivatives}
\label{subsec:SOSD}
We start by giving the definition of
the \emph{strong} second subderivatives of $J$.
\begin{definition}[strong second subderivative]
	\label{def:strong_star_subderivative}
	Let $w \in C_0(\Omega)$ and $x \in \MM(\Omega)$ be given.
	Then the strong second subderivative
	$J''(w, x; \cdot) \colon C_0(\Omega) \to [-\infty, \infty]$
	of $J$ at $w$ for $x$ is defined by
	\begin{equation*}
		J''(w, x; z)
		:=
		\inf
		\set*{
			\liminf_{k \to \infty} \frac{J(w + t_k z_k) - J(w) - t_k \dual{x}{z_k}}{t_k^2/2}
			\given
			t_k \searrow 0,
			z_k \to z
		}
		.
	\end{equation*}
	In the case that $J$ is Gâteaux differentiable at $w$,
	we will always use $x = J'(w)$ and, consequently,
	omit the argument $x$ of $J''$,
	i.e., we write $J''(w; \cdot)$ instead of $J''(w, x; \cdot)$.
	If for every $z \in C_0(\Omega)$ and every sequence $t_k \searrow 0$
	there exists a sequence $z_k \to z$ with
	\begin{equation*}
		J''(w, x; z)
		=
		\lim_{k \to \infty} \frac{J(w + t_k z_k) - J(w) - t_k \dual{x}{z_k}}{t_k^2/2}
		,
	\end{equation*}
	we say that $J$ is strongly-strongly twice epi-differentiable
	at $w$ for $x$.
\end{definition}
In order to obtain results for a large class of functions,
we investigate two prototypical situations:
\begin{itemize}
	\item
		$j(w) = \max(0, w)$
		with a kink at $w = 0$,
	\item
		$j \in C^{1,1}(\R)$ with directional derivatives of second order.
\end{itemize}

We start with the function $j = \max(0, \cdot)$.
In order to motivate the upcoming result,
we give a heuristic derivation based on the coarea formula
\[
 \int_\Omega \psi(x) \abs{\nabla u(x)} \d\lambda(x)  = \int_{\R} \int_{u^{-1}(s)} \psi(x)\d\HH^{d-1}(x) \d s,
\]
which holds for Lipschitz continuous $u$ and integrable $\psi$.
Here, $\HH^{d-1}$ is the usual $(d-1)$-dimensional Hausdorff measure.
Let $w \in C^1(\Omega) \cap C_0(\Omega)$
be given
such that $\set{w = 0}$ is a compact subset of $\Omega$
and assume that $\nabla w(\omega) \ne 0$ holds for all $\omega \in \set{w = 0}$.
This implies that $\set{w = 0}$ is a
$(d-1)$-dimensional
$C^1$-manifold.
Note that $\HH^{d-1}$ coincides with the surface measure on the manifold $\set{w = 0}$.
Since the Lebesgue measure of $\set{w = 0}$ is zero,
we get
$\partial J(w) = \set{x}$
with $x = \chi_{\set{w > 0}}$, see \cref{thm:rockit}.
For a second function $z \in C^1(\Omega) \cap C_c(\Omega)$,
we are going to calculate
\begin{align*}
	I &:=
	\lim_{t \searrow 0} \frac{J(w + t z) - J(w) - t \dual{x}{z}}{t^2 / 2}
	=
	\lim_{t \searrow 0}
	\int_\Omega
	\frac{ \max(0, w + t z ) - \max(0, w) - t \chi_{\set{w > 0}} z }{t^2 / 2}
	\d\lambda
	,
\end{align*}
cf.\ \cref{def:weak_star_subderivative}.
It is easy to check that the integrand vanishes
on the set $\set{ \sign(w) = \sign(w + t z)}$.
On its complement
$\Omega_t := \set{ \sign(w) \ne \sign(w + t z)}$,
the integrand is equal to $2 \abs{w + t z} / t^2$.
Now, we formally apply the coarea formula
with the functions $u = -w/z$,
$\psi = \abs{w + t z} / \abs{\nabla u}$
and obtain
\begin{equation*}
	I
	=
	\lim_{t \searrow 0} \int_{\Omega_t} \frac{2 \abs{w + t z }}{t^2} \d\lambda
	=
	\lim_{t \searrow 0}
	\int_0^t \int_{\set{u = \varepsilon}}
	\frac{2 \abs{w + t z }}{t^2 \abs{\nabla u}}
	\d\HH^{d-1} \d\varepsilon
	,
\end{equation*}
where we have used that $x\in \Omega_t$ if and only if $0<u(x)<t$.
Now, we plug in the values of $u$ and $\nabla u = -\nabla w / z + w \nabla z / z^2$,
and use $w = -\varepsilon z$ in the inner integral
to obtain
\begin{equation*}
	I
	=
	\lim_{t \searrow 0}
	\frac1t
	\int_0^t \int_{\set{w + \varepsilon z = 0}}
	\frac{2 (1 - \varepsilon/t) z^2}{\abs{\nabla w + \varepsilon \nabla z}}
	\d\HH^{d-1} \d\varepsilon
	=
	\lim_{t \searrow 0} \frac1t \int_0^t \parens*{1 - \frac{\varepsilon}{t}} \xi(\varepsilon) \d\varepsilon
	,
\end{equation*}
with
\begin{equation*}
	\xi(\varepsilon)
	=
	\int_{\set{w + \varepsilon z = 0}}
	\frac{2 z^2}{\abs{\nabla w + \varepsilon \nabla z}}
	\d\HH^{d-1}
	.
\end{equation*}
Since we integrate over a perturbation of the manifold $\set{w = 0}$,
it is reasonable to expect the continuity $\xi(\varepsilon) \to \xi(0)$ as $\varepsilon \searrow 0$.
Consequently, we expect to find
\begin{equation*}
	I
	=
	\lim_{t \searrow 0} \frac1t\int_0^t \parens*{1 - \frac{\varepsilon}{t}} \xi(\varepsilon) \d\varepsilon
	=
	\frac12 \xi(0)
	=
	\int_{\set{w = 0}} \frac{z^2}{\abs{\nabla w}} \d\HH^{d-1}
	.
\end{equation*}
	An anonymous reviewer has pointed out that this limit can be verified
	easily under the additional regularity $w,z \in C^2(\Omega)$
	by using the divergence theorem.
Our first goal is to derive this equation rigorously
with the relaxed regularity requirement $z \in C_0(\Omega)$,
see \cref{lemma:noncompstantin} below.

The upcoming theorem is prepared by the next lemma,
which studies the one-dimen\-sional case $d=1$.
\begin{lemma}
	\label{lem:diff_int_1d}
	Let $w \in C^1(\R)$ be given such that $\set{w = 0} = \set{0}$ and $w'(0) \ne 0$.
	Further, for $z \in C_c(\R)$ and $t > 0$ we define the set
	\begin{equation*}
		\Omega_t := \set{\sign(w) \ne \sign(w + t z)}.
	\end{equation*}
	\begin{enumerate}
		\item
			\label{itemlem:diff_int_1d_a}
			For each $\psi \in C(\R)$ we have for $t \searrow 0$
			\begin{subequations}
				\label{eq:diff_int_1d}
				\begin{align}
				\label{eq:diff_int_1d_a}
				\frac1t \int_{\Omega_t} \psi \d\lambda^1 &\to \frac{\psi(0) \abs{z(0)}}{\abs{w'(0)}}
				,
				\\
				\label{eq:diff_int_1d_b}
				\frac1t \int_{\Omega_t} \psi \sign(z) \d\lambda^1 &\to \frac{\psi(0) z(0)}{\abs{w'(0)}}
				,
				\\
				\label{eq:diff_int_1d_c}
				\frac1t \int_{\Omega_t} \frac{\abs{w}}{t} \d\lambda^1 &\to \frac12 \frac{z(0)^2}{\abs{w'(0)}}
				.
				\end{align}
			\end{subequations}
		\item
			\label{itemlem:diff_int_1d_b}
			If $\varepsilon, \delta, t_0 > 0$ are given such that
			$\varepsilon < \abs{w'(0)}$
			and such that
			\begin{subequations}
				\label{eq:conditions_eps_delta}
				\begin{align}
					\label{eq:conditions_eps_delta_1}
					\Omega_t &\subset [-\delta, \delta] \qquad\forall t \in (0,t_0)
					\\
					\label{eq:conditions_eps_delta_2}
					\abs{z(x) - z(0)} &\le \mrep{\varepsilon}{[-\delta,\delta]}  \qquad\forall \abs{x}\le\delta
					\\
					\label{eq:conditions_eps_delta_3}
					\abs{w(x) - w(0) - w'(0) x} &\le \mrep{\varepsilon \abs{x}}{[-\delta,\delta]} \qquad\forall \abs{x}\le\delta
				\end{align}
			\end{subequations}
			hold,
			then
		for all $t \in (0,t_0)$ we have
			\begin{subequations}
				\label{eq:bound_int_1d}
				\begin{align}
					\frac{\lambda^1(\Omega_t)}{t} &\le
					\frac{\abs{z(0)} + 2 \varepsilon}{\abs{w'(0)} - \varepsilon}
					,
					\label{eq:bound_int_1d_a}
					\\
					\frac1t \int_{\Omega_t} \frac{\abs{w}}{t} \d\lambda^1 &\le
					\frac{\parens{\abs{z(0)} + 2 \varepsilon}^2 (\abs{w'(0)} + \varepsilon)}{\parens{\abs{w'(0)} - \varepsilon}^2}
					.
					\label{eq:bound_int_1d_b}
				\end{align}
			\end{subequations}
	\end{enumerate}
\end{lemma}
\begin{proof}
	W.l.o.g.\ we consider the case $w'(0) > 0$,
	which implies $\sign(w(x))=\sign(x)$ for all $x\in \R$.
	First, we will prove \ref{itemlem:diff_int_1d_b}.
	Let $\varepsilon, \delta, t_0 > 0$ satisfying the assumptions in
	\ref{itemlem:diff_int_1d_b}
	be given.
	Further, let $t \in (0,t_0)$ be arbitrary.
	It is easy to check that
	\begin{subequations}
		\label{eq:diff_implications}
		\begin{align}
			\label{eq:diff_implications_1}
			0 &< x \le t \frac{-z(0) - \varepsilon}{w'(0) + \varepsilon}
			&\Rightarrow\quad
			x &\in \Omega_t \cap (0,\infty)
			&\Rightarrow\quad
			0 &< x \le t \frac{-z(0) + \varepsilon}{w'(0) - \varepsilon}
			,
			\\
			\label{eq:diff_implications_2}
			0 &> x \ge t \frac{-z(0) + \varepsilon}{w'(0) + \varepsilon}
			&\Rightarrow\quad
			x &\in \Omega_t \cap (-\infty, 0)
			&\Rightarrow\quad
			0 &> x \ge t \frac{-z(0) - \varepsilon}{w'(0) - \varepsilon}
			.
		\end{align}
	\end{subequations}
	Indeed, the second implication follows from the chain of inequalities
	\begin{equation*}
		0
		\ge
		w(x) + t z(t)
		\ge
		w(0) + w'(0) x - \varepsilon x + t z(0) - t \varepsilon
		=
		(w'(0) - \varepsilon) x + t \, (z(0) - \varepsilon)
	\end{equation*}
	for $x \in \Omega_t \cap (0,\infty)$.
	The other implications follow similarly.
	From \eqref{eq:diff_implications_1} and \eqref{eq:diff_implications_2}
	we get
	\begin{equation*}
		\frac{\lambda^1(\Omega_t)}{t}
		\le
		\frac{\max\set{0, \varepsilon - z(0)} + \max\set{0, z(0) + \varepsilon}}{w'(0) - \varepsilon}
		\le
		\frac{\abs{z(0)} + 2 \varepsilon}{w'(0) - \varepsilon},
	\end{equation*}
	i.e., \eqref{eq:bound_int_1d_a}.
	From \eqref{eq:diff_implications_1} and \eqref{eq:diff_implications_2}
	we have
	$
		\abs{w(x)}
		\le
		(w'(0) + \varepsilon) \abs{x}
		\le
		(w'(0) + \varepsilon) t \frac{\abs{z(0)} + \varepsilon}{w'(0) - \varepsilon}.
	$
	Together with \eqref{eq:bound_int_1d_a},
	this shows
	\ref{itemlem:diff_int_1d_b}.

	In order to validate \ref{itemlem:diff_int_1d_a},
	let $\varepsilon \in (0,w'(0))$ be arbitrary.
	Then, there exists $\delta > 0$, such that
	\eqref{eq:conditions_eps_delta_2} and \eqref{eq:conditions_eps_delta_3}
	are satisfied.
	Let $K$ be the support of $z$ and define
	\begin{align*}
		c &:= \inf\set{ \abs{w(x)} \given x \in K \setminus (-\delta, \delta) } > 0,
		&
		C &:= \sup\set{ \abs{z(x)} \given x \in K } < \infty.
	\end{align*}
	Then, for $t_0 := c / C$ we have
	$\sign(w(x)) =  \sign(w(x) + t z(x))$ for all $x \in K \setminus (-\delta, \delta)$
	and this shows
	\eqref{eq:conditions_eps_delta_1}.
	Thus, the first part of the proof shows \eqref{eq:diff_implications}.
	From the continuity of $\psi$, \eqref{eq:diff_int_1d_a} and \eqref{eq:diff_int_1d_b}
	follow easily.
	It remains to show \eqref{eq:diff_int_1d_c}.
	In case $z(0) = 0$, this follows directly from \eqref{eq:bound_int_1d_b}.
	We focus on the case $z(0) > 0$, the remaining case $z(0) < 0$
	is similar.
	W.l.o.g., we assume $\varepsilon < z(0)$.
	Since the integrand is positive, \eqref{eq:diff_implications} implies
	\begin{align*}
		\frac12 (w'(0) - \varepsilon) \parens*{t \frac{-z(0)-\varepsilon}{w'(0)+\varepsilon}}^2
		&=
		\int_0^{t \frac{-z(0)-\varepsilon}{w'(0)+\varepsilon}}
		(w'(0) - \varepsilon) x
		\d\lambda^1(x)
		\le
		\int_{\Omega_t} \abs{w} \d\lambda^1
		\\&
		\le
		\int_0^{t \frac{-z(0)+\varepsilon}{w'(0)-\varepsilon}}
		(w'(0) + \varepsilon) x
		\d\lambda^1(x)
		=
		\frac12 (w'(0) + \varepsilon) \parens*{t \frac{-z(0)+\varepsilon}{w'(0)-\varepsilon}}^2
		.
	\end{align*}
	After division by $t^2$,
	the left-hand side and the right-hand side converge
	to $z(0)^2/(2 w'(0))$ as $\varepsilon \searrow 0$.
	This shows \eqref{eq:diff_int_1d_c}.
\end{proof}

Using standard coordinate transform arguments,
this result can be lifted to the $d$-dimensional situation.

\begin{theorem}
	\label{thm:taylor_expansion_integral_rn}
	Let $w \in C^1(\Omega)$ be given such that
	$\nabla w \ne 0$ on $\set{w = 0}$.
	Further, let $z \in C_c(\Omega)$ be given.
	For all $t > 0$, we define $\Omega_t := \set{ \sign(w) \ne \sign(w + t z)}$.
	Then, for all $\psi \in C(\Omega)$ we have
	\begin{subequations}
		\label{eq:diff_int_rn}
		\begin{align}
			\label{eq:diff_int_rn_a}
			\frac1t \int_{\Omega_t} \psi \d\lambda &\to \int_{\set{w=0}}\frac{\psi \abs{z}}{\abs{\nabla w}} \d\HH^{d-1}
			,
			\\
			\label{eq:diff_int_rn_b}
			\frac1t \int_{\Omega_t} \psi \sign(z) \d\lambda &\to \int_{\set{w=0}}\frac{\psi z}{\abs{\nabla w}} \d\HH^{d-1}
			,
			\\
			\label{eq:diff_int_rn_c}
			\frac1t \int_{\Omega_t} \frac{\abs{w}}{t} \d\lambda &\to \frac12 \int_{\set{w=0}}\frac{z^2}{\abs{\nabla w}} \d\HH^{d-1}
		\end{align}
	\end{subequations}
	as $t \searrow 0$.
\end{theorem}
\begin{proof}
	We prove \eqref{eq:diff_int_rn_a},
	the remaining limits can be verified by analogous arguments.
	We start with a local result.
	We set $Z = \set{w = 0} \cap \supp z$, which is a compact subset of $\Omega$.
	Let a point $p \in Z$,
	a bounded open set $B \subset \R^{d-1}$, an open interval $I = (a,b)$
	and $\gamma \in C^1(\bar B)$ be given
	such that
	\begin{subequations}
	\label{eq_thm32_local}
	\begin{align}
		p &\in  B\times I
		,\qquad
		\overline{B \times I} \subseteq \Omega
		, \qquad
		Z \cap   (B \times I)  =  \set{(x,\gamma(x)) \given x \in B }
		\\
		W  &:= \set{(x,y) \given x \in B , \abs{y - \gamma (x)} < \varepsilon } \subseteq B  \times I
		, \qquad
		\partial_d w(p) > 0
	\end{align}
	\end{subequations}
	for some $\varepsilon > 0$.
	Further, let $\varphi \in C_c(B \times I)$ be given.

	For all $x \in B$ we define the cross section
	\begin{equation*}
		\Omega_{t,x}
		:=
		\set{ y \in I \given (x,y) \in \Omega_t}
		=
		\set{ y \in I \given \sign(w(x,y)) \ne \sign( (w+t z)(x,y))}
		.
	\end{equation*}
	By Fubini, we have
	\begin{equation*}
		\frac1t \int_{\Omega_t \cap (B \times I)} \psi \varphi \d\lambda
		=
		\int_B \frac1t \int_{\Omega_{t,x}} \psi(x,y) \varphi(x,y) \d\lambda^1(y) \d\lambda^{d-1}(x)
		.
	\end{equation*}
	Now, we can apply \eqref{eq:diff_int_1d_a}
	for all $x \in B$ and obtain
	\begin{equation*}
		\frac1t \int_{\Omega_{t,x}} \psi(x,y) \varphi(x,y) \d\lambda^1(y)
		\to
		\frac{\psi(x,\gamma(x)) \varphi(x,\gamma(x)) \abs{z(x,\gamma(x))}}{\abs{\partial_d w(x,\gamma(x))} }
		\qquad\forall x \in B.
	\end{equation*}
	In order to pass to the limit in the outer integral,
	we need an upper bound.
	To this end,
	let $0 < \varepsilon < \inf\set{\abs{\partial w_d (x,\gamma(x))} \given x \in B}$
	be given.
	By uniform continuity, we find $\delta > 0$,
	such that
	\begin{align*}
		\abs{z(x, \gamma(x) + h) - z(x,\gamma(x))} &\le \varepsilon \phantom{\abs{h}} \qquad\forall \abs{h} \le \delta, x \in B
		\\
		\abs{w(x, \gamma(x) + h) - w(x,\gamma(x)) - \partial_d w(x,\gamma(x)) h} &\le \varepsilon \abs{h} \qquad\forall \abs{h} \le \delta, x \in B
	\end{align*}
	hold.
	Similarly to the proof of \cref{lem:diff_int_1d}~\ref{itemlem:diff_int_1d_a},
	we find $t_0 > 0$, such that
	$\Omega_{t,x} \subset \gamma(x) + [-\delta,\delta]$
	for all $t\in (0,t_0)$.
	Thus, we can apply \cref{lem:diff_int_1d}~\ref{itemlem:diff_int_1d_b},
	and this yields
	the integrable bound
	\begin{equation*}
		\frac1t \int_{\Omega_{t,x}} \psi(x,y) \varphi(x,y) \d\lambda^1(y)
		\le
		\frac1t \lambda^1(\Omega_{t,x}) M
		\le
		M \frac{\abs{z(x,\gamma(x))} + 2 \varepsilon}{\abs{\partial_d w(x,\gamma(x))} - \varepsilon}
	\end{equation*}
	where $M$ is an upper bound for $\abs{\psi \varphi}$.
	Thus, we can apply the dominated convergence theorem to obtain
	\begin{equation*}
		\frac1t \int_{\Omega_t \cap (B \times I)} \psi \varphi \d\lambda
		\to
		\int_B
		\frac{\psi(x,\gamma(x)) \varphi(x,\gamma(x)) \abs{z(x,\gamma(x))}}{\abs{\partial_d w(x,\gamma(x))} }
		\d\lambda^{d-1}(x)
		.
	\end{equation*}
	By differentiating $w(x,\gamma(x)) = 0$,
	we obtain
	$-\partial_i w(x,\gamma(x)) = \partial_d w(x,\gamma(x)) \partial_i \gamma(x)$
	for all $i = 1,\ldots, d-1$.
	Thus,
	\begin{equation*}
		\abs{\nabla w(x,\gamma(x))}^2
		=
		\abs{\partial_d w(x,\gamma(x))}^2 (\abs{\gamma'(x)}^2 + 1).
	\end{equation*}
	This yields
	\begin{align*}
		\int_B
		\parens*{\frac{\psi \varphi \abs{z}}{\abs{\partial_d w} }}|_{(x,\gamma(x))}
		\d\lambda^{d-1}(x)
		&=
		\int_B
		\parens*{\frac{\psi \varphi \abs{z}}{\abs{\nabla w} }}|_{(x,\gamma(x))}
		\sqrt{1 + \abs{\gamma'(x)}^2}
		\d\lambda^{d-1}(x)
		\\&=
		\int_{Z \cap (B \times I)}
		\frac{\psi \varphi \abs{z}}{\abs{\nabla w} }
		\d\HH^{d-1}
		.
	\end{align*}
The compact set $Z$ can be covered by finitely many, possibly rotated and translated sets $B\times I$ of type \eqref{eq_thm32_local}.
The proof is then finished using a standard partition-of-unity argument.
\end{proof}
As in \cref{thm:taylor_expansion_integral_rn},
we will often work with
$w \in C^1(\Omega)$
such that $\set{ w = 0 } \cap \set{ \nabla w = 0} = \emptyset$.
This condition means $\nabla w \ne 0$ everywhere on $\set{ w = 0 }$.
From Stampacchia's lemma, we have $\nabla w = 0$ a.e.\ on $\set{w = 0}$.
Thus, $\lambda(\set{ w = 0 }) = 0$.
We will now apply \cref{thm:taylor_expansion_integral_rn} to obtain a second-order expansion of $w\mapsto\max(w,0)$.

\begin{lemma}
\label{lemma:constantin}
Define $j\colon\R\to\R$ by $j(w):=\max(w,0)$.
 Let a function $w\in C^1(\Omega)$ be given such that $\set{w=0}\cap \set{\nabla w=0}=\emptyset$.
 Let $z_k \to z$ in $C(\bar \Omega)$ with $\supp z_k \subset K$ for some compact $K\subset\Omega$.
 Let $t_k \searrow0$.
 Then we have
 \begin{equation} \label{eq_limit_z_compact_support}
 \lim_{k\to\infty}
 \frac2{t_k^2} \int_\Omega j(w + t_kz_k) - j(w) - t_k j'(w;z_k) \d\lambda
 =
 \int_{\{w=0\}} \frac{ z^2 }{\abs{\nabla w}} \d\HH^{d-1}
 \in [0,\infty)
 \end{equation}
for $k\to\infty$.
\end{lemma}
\begin{proof}
	We first consider the situation that $z_k$ does not depend on $k$.
	For every $k \in \N$, we define the set
	$\Omega_k := \set{ \sign(w) \ne \sign(w + t_k z)}$
	and have
	\begin{equation*}
		\int_\Omega j(w + t_k z) - j(w) - t_k j'(w; z) \d\lambda
		=
		\int_{\Omega_k} \abs{w + t_k z} \d\lambda
		=
		\int_{\Omega_k} t_k \abs{z} - \abs{w} \d\lambda
		.
	\end{equation*}
	Now, the claim follows from \eqref{eq:diff_int_rn_b} and \eqref{eq:diff_int_rn_c}.

	For the general case, we consider
	\begin{align*}
		I_k &:=
		\int_\Omega j(w + t_k z_k) - j(w) - t_k j'(w; z_k) \d\lambda
		-
		\int_\Omega j(w + t_k z) - j(w) - t_k j'(w; z) \d\lambda
		\\&
		=
		\int_\Omega j(w + t_k z_k) - j(w + t_k z) - t_k j'(w; z_k - z) \d\lambda
		.
	\end{align*}
	Here, we used $\lambda(\set{ w = 0 }) = 0$, which implies $j'(w; z_k)-j'(w; z)=j'(w; z_k - z)$.
	The absolute value of the integrand can be estimated
	by $2 t_k \abs{z_k - z}$,
	and it vanishes on the complement of
	$A_k := \set{ \abs{w} \le (\norm{z}_{L^\infty(\lambda)} + \norm{z_k}_{L^\infty(\lambda)}) t_k } \cap K$.
	By \cite[Lemma~3.2]{DeckelnickHinze2012}, there exists $c>0$ such that
	$\lambda(A_k) \le c (\norm{z}_{L^\infty(\lambda)} + \norm{z_k}_{L^\infty(\lambda)}) t_k$.
	Hence, it follows
	\[
	 \frac2{t_k^2} \abs{I_k}  \le \frac{4}{t_k}\lambda(A_k)  \norm{z_k - z}_{L^\infty(\lambda)} \to 0,
	\]
	which proves the claim for the general case.
\end{proof}

In the next step, we want to drop the assumption that $z$ is compactly supported in $\Omega$.
Interestingly, the result of the previous \cref{lemma:constantin} is not valid
anymore in this case. In particular, we can no longer choose $z_k\equiv z$ in \eqref{eq_limit_z_compact_support},
as the following example shows.

\begin{example}
	\label{ex:no_convergence_zk_z}
	We take $\Omega = (0,2) \subset \R$.
	On $(0,1)$, the functions $w$ and $z$
	are given by $w(s) = s^3$ and $z(s) = -s$.
	We extend these functions to $(0,2)$, such that
	$w \in C^1(\Omega) \cap C_0(\Omega)$, $z \in C_0(\Omega)$
	and $w > 0$, $w + z > 0$ on $(1,2)$.
	For $t \in (0,1)$, we set
	\begin{equation*}
		\Omega_t
		:=
		\set{\sign(w) \ne \sign(w + t z)}
		=
		(0,\sqrt{t}].
	\end{equation*}
	A short computation shows
	\begin{equation*}
		\int_\Omega j(w + tz) - j(w) - t j'(w;z) \d\lambda
		=
		\int_0^{\sqrt{t}}
		t s - s^3
		\d\lambda(s)
		=
		\frac{t^2}{4}.
	\end{equation*}
	Hence
	\begin{equation*}
		\lim_{t \to 0} \frac2{t^2} \int_\Omega j(w + tz) - j(w) - t j'(w;z) \d\lambda = \frac12.
	\end{equation*}
	However, the boundary integral $\int_{\set{w = 0}} \frac{ z^2 }{\abs{\nabla w}} \d\HH^{d-1}$
	vanishes due to $\set{w = 0} = \emptyset$.
	Hence, \eqref{eq_limit_z_compact_support} is not true for $z_k\equiv z$.
\end{example}

This example shows that we cannot expect \eqref{eq_limit_z_compact_support} to hold for all sequences $(z_k)$ converging
to some non-compactly supported $z$.
In the next lemma, we show that there exists a sequence $(z_k) \subset C_c(\Omega)$ with $ z_k \to z$ in $C_0(\Omega)$
such that equality holds in the limit \eqref{eq_limit_z_compact_support}.
In the general case, we have to replace in \eqref{eq_limit_z_compact_support} the limit by the limit inferior and the equality sign by an inequality sign.

\begin{lemma}
	\label{lemma:noncompstantin}
	Define $j\colon\R\to\R$ by $j(w):=\max(w,0)$.
	Let a function $w\in C^1(\Omega)$ be given such that $\{w=0\}\cap \{\nabla w=0\}=\emptyset$.
	Let $z \in C_0(\Omega)$ be given.
	For all sequences $t_k \searrow 0$ and $z_k \to z$ in $C_0(\Omega)$,
	we have
	\begin{equation*}
		\liminf_{k \to \infty} \frac2{t_k^2} \int_\Omega j(w + t_kz_k) - j(w) - t_k j'(w;z_k) \d\lambda
		\ge
		\int_{\{w=0\}} \frac{ z^2 }{\abs{\nabla w}} \d\HH^{d-1}
		\in
		[0,\infty].
	\end{equation*}
	Moreover, for every sequence $t_k \searrow 0$ there exists $C_c(\Omega) \ni z_k \to z$ in $C_0(\Omega)$,
	such that
	\begin{equation*}
		\lim_{k \to \infty} \frac2{t_k^2} \int_\Omega j(w + t_kz_k) - j(w) - t_k j'(w;z_k) \d\lambda
		=
		\int_{\{w=0\}} \frac{ z^2 }{\abs{\nabla w}} \d\HH^{d-1}
		\in
		[0,\infty].
	\end{equation*}
\end{lemma}
\begin{proof}
	Let $t_k \searrow 0$ and $z_k \to z$ in $C_0(\Omega)$ be arbitrary.
	For every $\psi \in C_c(\Omega)$, $0 \le \psi \le 1$,
	we have
	\begin{align*}
		\int_\Omega j(w + t_kz_k) - j(w) - t_k j'(w;z_k) \d\lambda
		& \ge
		\int_\Omega \psi ( j(w + t_k  z_k) - j(w) - t_k j'(w;z_k)) \d\lambda \\
		& \ge
		\int_\Omega j(w + t_k \psi z_k) - j(w) - t_k j'(w;\psi z_k) \d\lambda
		,
	\end{align*}
	where we used the non-negativity of the integrand in the first inequality, and the convexity of $j$ and positive homogeneity of $j'$ in
	the second inequality.
	For the term on the right-hand side, we can apply \cref{lemma:constantin}
	and obtain
	\begin{equation*}
		\liminf_{k \to \infty} \frac2{t_k^2} \int_\Omega j(w + t_kz_k) - j(w) - t_k j'(w;z_k) \d\lambda
		\ge
		\int_{\{w=0\}} \frac{ (\psi z)^2 }{\abs{\nabla w}} \d\HH^{d-1}
		.
	\end{equation*}
	Taking the supremum over all possible $\psi$ yields the first claim.

	To address the second claim, we start by a sequence $C_c(\Omega) \ni \hat z_n \to z$
	such that $\abs{\hat z_n}$ is monotonically increasing.
	We use a diagonal sequence argument as in \cite[Lemma~2.12(ii)]{ChristofWachsmuth2019}:
	Due to \cref{lemma:constantin},
	we find a strictly increasing sequence $K_n$
	such that
	\begin{equation*}
		\abs*{
			\frac2{t_k^2} \int_\Omega j(w + t_k \hat z_n) - j(w) - t_k j'(w; \hat z_n) \d\lambda
			-
			\int_{\{w=0\}} \frac{ \hat z_n^2 }{\abs{\nabla w}} \d\HH^{d-1}
		}
		\le
		\frac1n
	\end{equation*}
	holds for all $n \in \N$ and $k \ge K_n$.
	We redefine $K_1 := 1$
	and set $n_k := \sup\set{ n \in \N \given   K_{n} \le k }$ for all $k \in \N$,
	which is finite, since $K_n\to\infty$.
	Then, it holds $k \geq K_{n_k}$ for all $k$ by definition
	and $n_k \to \infty$ monotonously for $k \to \infty$.
	Now, we define $z_k := \hat z_{n_k}$
	and this yields
	\begin{equation*}
		\abs*{
			\frac2{t_k^2} \int_\Omega j(w + t_k z_k) - j(w) - t_k j'(w; z_k) \d\lambda
			-
			\int_{\{w=0\}} \frac{ z_k^2 }{\abs{\nabla w}} \d\HH^{d-1}
		}
		\le
		\frac1{n_k}
		\to
		0,
	\end{equation*}
	where the inequality holds for all $k \in \N$ with $n_k \ge 2$.
	Due to $n_k \to \infty$, we have $z_k \to z$ in $C_0(\Omega)$.
	Since $(n_k)$ is monotonically increasing,  $\abs{z_k}$ is monotonically increasing as well.
	Thus, the monotone convergence theorem
	yields
	\begin{equation*}
		\int_{\set{w=0}} \frac{ z_k^2 }{\abs{\nabla w}} \d\HH^{d-1}
		\to
		\int_{\set{w=0}} \frac{ z  ^2 }{\abs{\nabla w}} \d\HH^{d-1}
		.
	\end{equation*}
	This proves the second claim.
\end{proof}

This lemma shows that $J$ is strongly-strongly twice epi-differentiable
in the sense of \cref{def:strong_star_subderivative}
and
the strong second subderivative is given by
\begin{equation}
	\label{eq:strong_second_order_subderivative}
	J''(w; z) =
	\int_{\{w=0\}} \frac{ z^2 }{\abs{\nabla w}} \d\HH^{d-1}
	.
\end{equation}
	A similar result has been obtained in
	\cite[Theorem~5.2.14]{Christof2018}
	without using the strong regularity condition
	$\set{w = 0} \cap \set{\nabla w = 0} = \emptyset$.
	Therein, the constructed recovery sequence
	converges in $H_0^1(\Omega)$,
	whereas \cref{lemma:noncompstantin}
	gives uniform convergence.

Next, we address a function $j$ without kinks.
\begin{lemma}\label{lemma:zweeteableitung}
 Let $j\in C^1(\R)$ with  locally Lipschitz continuous  $j'$ and such that
 \[
  j''(w;z):=\lim_{\substack{t\searrow0 \\ \hat z \to z}} \frac2{t^2} (j(w+t\hat z)- j(w) - t j'(w)\hat z) \in \R
 \]
 exists for all $w,z\in \R$.

 Let $w\in C^1(\Omega)$ be bounded.
 Let $z_k \to z$ in $C(\bar \Omega)$.
 Let $t_k \searrow0$.
 Then we have
 \[
 \frac2{t_k^2} \int_\Omega j(w + t_kz_k) - j(w) - t_k j'(w)z_k \d\lambda
 \to
 \int_\Omega j''(w;z) \d\lambda  \in \R
 \]
for $k\to\infty$.
\end{lemma}
\begin{proof}
By the boundedness of $w$ and the Lipschitz continuity of $j'$, $\frac2{t_k^2} (j(w+t_kz_k)- j(w) - t_k j'(w)z_k)$ is uniformly bounded.
Then the claim follows by dominated convergence.
\end{proof}

\begin{remark}
 It is an open question, whether the integral functional in \cref{lemma:zweeteableitung} is strongly-strongly
 twice
 epi-differentiable on $C_0(\Omega)$ if $j$ is assumed to be second-order epi-differentiable, only.
 In particular, the existence of a recovery sequence seems to be impossible to prove.
 For differentiation of integral functionals on $L^p$-spaces, this was done in \cite[Section~5]{Do} for convex and  in \cite{Levy1993} for non-convex functions.
 In order to construct the recovery sequence, a selection theorem for measurable multi-functions was used \cite[Thm. 1.4]{Levy1993}.
 It is unclear, under which assumptions on $j$ this multi-function allows for a continuous selection
 without assuming continuity of $w\mapsto j''(w;z)$.
\end{remark}

\begin{remark}
	Since $j \colon \R\to \R$, we can omit the limit ``$\hat z\to z$'' in the definition of $j''(w;z)$, i.e.,
	it is sufficient to require
 \[
  j''(w;z)=\lim_{t\searrow0} \frac2{t^2} (j(w+t z)- j(w) - t j'(w)z) \in \R.
 \]
 By construction, $z \mapsto j''(w;z)$ is positively $2$-homogeneous. Hence, $j''(w;\cdot)$ is determined by
 the values $j''(w;z)$ for $z\in \{-1,+1\}$.
\end{remark}

By combining the previous two lemmas,
we obtain the convergence of a second-order difference quotient associated with piecewise smooth functions.
\begin{assumption}
 \label{asm:dzweij_pw}
 Let $j_0\in C^1(\R)$ with  locally Lipschitz continuous  $j_0'$ be given such that
 \[
  j_0''(w;z):=\lim_{t\searrow0} \frac2{t^2} (j_0(w+t z)- j_0(w) - t j_0'(w) z) \in \R
 \]
 exists for all $w,z\in \R$.
 Let $m\in \N$, positive numbers $a_i$, $i=1\dots m$,  and distinct numbers $b_i\in \R$, $i=1\dots m$ be given.
 We define $j \colon \R  \to \R$ and $J \colon C_0(\Omega) \to \R$
 via
 \begin{align*}
  j(w)
  &:=
  j_0(w)
  +
  \sum_{i = 1}^m a_i \abs{ w - b_i }
  ,
  &
  J(w)
  &:=
  \int_\Omega
  j(w(\omega))
  \d \lambda(\omega)
  .
 \end{align*}
\end{assumption}

\begin{theorem}\label{thm_dzweij_pw}
 We consider the setting of \cref{asm:dzweij_pw}.
Let a bounded $w\in C^1(\Omega)$ be given such that
\begin{equation}\label{eq_ass_w_pw}
 \{w=b_i \}  \cap \{\nabla w=0\}=\emptyset \qquad \forall i=1\dots m.
\end{equation}
Then $J$ is (Gâteaux) differentiable at $w$ with
\[
 J'(w)z = \int_\Omega \parens[\bigg]{ j_0'(w) + \sum_{i = 1}^m a_i \sign( w - b_i ) } z  \d\lambda.
\]
For $z\in C_0(\Omega)$, we define
\begin{equation}
\label{eq_dzweij_pw}
 J''(w; z):=
 \int_\Omega j_0''(w;z) \d\lambda
 +
 \sum_{i=1}^m
		2a_i \int_{\{w=b_i\}} \frac{ z^2 }{\abs{\nabla w}} \d\HH^{d-1}
		\in
		(-\infty,\infty].
\end{equation}
	Let $z \in C_0(\Omega)$ be given.
	For all sequences $t_k \searrow 0$ and $z_k \to z$ in $C_0(\Omega)$,
	we have
	\begin{equation*}
		\liminf_{k \to \infty} \frac{J(w + t_kz_k) - J(w) - t_k J'(w)z_k}{t_k^2 / 2}
		\ge  J''(w; z).
	\end{equation*}
	Moreover, for every sequence $t_k \searrow 0$ there exists $C_c(\Omega) \ni z_k \to z$ in $C_0(\Omega)$,
	such that
	\begin{equation*}
		\lim_{k \to \infty} \frac{J(w + t_kz_k) - J(w) - t_k J'(w)z_k}{t_k^2 / 2}
		=  J''(w; z).
	\end{equation*}
	Thus, $J$ is strongly-strongly twice epi-differentiable at $w$ for $J'(w)$ with
\(
J''(w,J'(w); z):= J''(w; z)
\).
\end{theorem}
\begin{proof}
We can write $\abs{ w - b_i } = 2\max( w - b_i,0) - ( w - b_i)$. The functions $w_i:=w-b_i$ satisfy
the requirements of  \cref{lemma:noncompstantin}.
The first claim is now a direct consequence of \cref{lemma:noncompstantin} and \cref{lemma:zweeteableitung}.
To prove the second claim, i.e., the existence of a recovery sequence, we use the diagonal sequence argument of the proof of \cref{lemma:noncompstantin}.
Let  $C_c(\Omega) \ni \hat z_n \to z$ be a sequence such that $\abs{\hat z_n}$ is monotonically increasing.
Define $j_i(w):=\max( w - b_i,0)$.
Then define a strictly increasing sequence $K_n\in \N$ with the property that
\begin{equation*}
	\abs*{
		\frac2{t_k^2} \int_\Omega j_i(w + t_k \hat z_n) - j_i(w) - t_k j_i'(w; \hat z_n) \d\lambda
		-
		\int_{\{w=b_i\}} \frac{ \hat z_n^2 }{\abs{\nabla w}} \d\HH^{d-1}
	}
	\le
	\frac1n
\end{equation*}
holds for all $i=1\dots m$, $n \in \N$, and $k \ge K_n$.
This is possible due to \cref{lemma:constantin}. Now we can conclude as in the proof of \cref{lemma:noncompstantin}.
Due to \cref{lemma:zweeteableitung} we get the convergence of the second-order difference quotients for the full functional $J$.
\end{proof}

\begin{corollary}
 We consider the setting of \cref{asm:dzweij_pw}
 with a convex $j_0$.
 Let a bounded $w\in C^1(\Omega)$ satisfy \eqref{eq_ass_w_pw}.
 Then the mapping $z\mapsto J''(w;z)$ is convex.
\end{corollary}

It remains to compute the convex conjugate of $z\mapsto J''(w;z)$.
For a measurable subset $\ZZ \subset \Omega$,
we utilize the measure $\HH^{d-1} |_\ZZ$ given by
\begin{equation*}
	\HH^{d-1}|_{\ZZ}(A):=\HH^{d-1}(A\cap \ZZ)
	\qquad
	\forall A \in \BB(\Omega).
\end{equation*}

\begin{theorem}\label{theo312}
 \label{thm:conj_subderivative}
 We consider the setting of \cref{asm:dzweij_pw}
 with a convex $j_0$.
 Let a bounded $w\in C^1(\Omega)$ satisfy \eqref{eq_ass_w_pw}.
 We set
 \begin{equation}
  \label{eq:def_ZZ}
  \ZZ:= \bigcup_{i=1}^m \set{w=b_i}
 \end{equation}
and define $a \colon \ZZ\to \R$ by
\[
 a(\omega) = a_i \text{ if } w(\omega) = b_i.
\]
Assume
\begin{equation}\label{eq_int_nablaw_z_finite}
 \int_\ZZ \frac1{\abs{\nabla w}}\d\HH^{d-1}<+\infty.
\end{equation}
Let $\mu \in \MM(\Omega)$.
If there exist densities $v_1 \in L^1(\lambda)$, $v_2 \in L^1(\HH^{d-1}|_\ZZ)$
	such that
	$\mu = v_1 \lambda + v_2 \HH^{d-1}|_{\ZZ}$
then
	\begin{equation*}
		\parens*{\frac12 J''(w;\cdot)}\conjugate(\mu)
		=
		\frac12  \int_{\ZZ} \frac{\abs{\nabla w}}{2a} v_2^2 \d\HH^{d-1}
		+
		\int_\Omega \left(\frac12 j_0''(w;\cdot) \right)\conjugate(v_1)  \d\lambda \in [0,+\infty].
	\end{equation*}
Otherwise  $\parens*{\frac12 J''(w;\cdot)}\conjugate(\mu) = + \infty$.
\end{theorem}
\begin{proof}
	We are going to apply \cref{thm:rockit}.
	To this end, we have to write $\frac12 J''(w;\cdot)$ as an integral functional.
	Therefore, we define the measure
	\begin{equation*}
		\nu := \lambda + \HH^{d-1}|_{\ZZ}
	\end{equation*}
	and the integrand $\iota \colon \Omega \times \R \to \R$ by
	\begin{equation*}
		\iota(\omega, z) := \frac12
		\begin{cases}
			\frac{2a(\omega)}{\abs{\nabla w(\omega)}} z^2
			& \text{if } \omega \in \ZZ,
			\\
			 j_0''(w(\omega);z)
			& \text{if } \omega \not\in \ZZ.
		\end{cases}
	\end{equation*}
	Note that we have $\lambda(\ZZ) = 0$.
	Thus
	\begin{align*}
		\int_\Omega \iota(\omega, z(\omega)) \d\nu(\omega)
		&=
		\int_\ZZ \iota(\omega, z(\omega)) \d\nu(\omega)
		+
		\int_{\Omega\setminus\ZZ} \iota(\omega, z(\omega)) \d\nu(\omega)
		\\&=
		\int_\ZZ \iota(\omega, z(\omega)) \d\HH^{d-1}|_{\ZZ}(\omega)
		+
		\int_{\Omega\setminus\ZZ} \iota(\omega, z(\omega)) \d\lambda(\omega)
		\\&=
		 \frac12\int_\ZZ \frac{2a(\omega)}{\abs{\nabla w(\omega)}} z^2 \d\HH^{d-1}|_{\ZZ}(\omega)
		+
		\frac12 \int_{\Omega\setminus\ZZ} j_0''(w(\omega), z(\omega)) \d\lambda(\omega)
		\\&=
		 \frac12\int_\ZZ \frac{2a(\omega)}{\abs{\nabla w(\omega)}} z^2 \d\HH^{d-1}(\omega)
		+
		\frac12 \int_{\Omega} j_0''(w(\omega), z(\omega)) \d\lambda(\omega)
		\\&=
		\frac12 J''(w;z)
		.
	\end{align*}
	It is clear that $\nu$ satisfies the assumptions of \cref{thm:rockit}
	and that $\iota$ is a convex Carathéodory function.
	The integrability $\iota(\cdot, x) \in L^1(\nu)$ for all $x \in \R$
	is a consequence of \eqref{eq_int_nablaw_z_finite}.
	Thus, all assumptions of \cref{thm:rockit} are satisfied
	and therefore, we obtain
	\begin{equation*}
		\parens*{\frac12 J(w;\cdot) }\conjugate(\mu)
		=
		\begin{cases}
		\int_\Omega \iota\conjugate\parens*{\omega , \frac{\d\mu}{\d\nu}(\omega)} \d\nu(\omega)
		&\text{if } \mu \ll \nu,
		\\
		+\infty & \text{otherwise}
		\end{cases}
	\end{equation*}
	for all $\mu \in \MM(\Omega)$.
	Next, we compute the conjugate of $\iota(\omega,\cdot)$.
	From the definition of $\iota$,
	we get
	\begin{equation*}
		\iota\conjugate(\omega,y) :=
		\begin{cases}
			\frac12  \frac{\abs{\nabla w(\omega)}}{2a(\omega)} y^2
			& \text{if } \omega \in \ZZ,
			\\
			\left(\frac12 j_0''(w(\omega);\cdot) \right)\conjugate(y)
			& \text{if } \omega \not\in \ZZ.
		\end{cases}
	\end{equation*}
	Since the measures $\lambda$ and $\HH^{d-1}|_{\ZZ}$
	are singular,
	it is clear that $\mu \ll \nu$
	if and only if
	\begin{equation*}
		\mu|_{\Omega \setminus \ZZ} \ll \lambda
		\quad\text{and}\quad
		\mu|_{\ZZ} \ll \HH^{d-1}|_{\ZZ}
		.
	\end{equation*}
	This is, in turn, equivalent
	to the existence of densities
	$v_1 \in L^1(\lambda)$, $v_2 \in L^1(\HH^{d-1}|_{\ZZ})$
	such that
	$\mu = v_1 \lambda + v_2 \HH^{d-1}|_{\ZZ}$.
	In this case we have
	\begin{equation*}
		\parens*{\frac12 J(w;\cdot) }\conjugate(\mu)
		=
		\int_\Omega \iota\conjugate\parens*{\omega , v_1(\omega)} \d\lambda(\omega)
		+
		\int_\ZZ \iota\conjugate\parens*{\omega , v_2(\omega)} \d\HH^{d-1}(\omega)
		.
	\end{equation*}
	From the structure of $\iota\conjugate$,
	the assertion of the structure of $\parens*{\frac12 J(w;\cdot) }\conjugate(\mu)$ follows.

Since $\left(\frac12 j_0''(w(\omega);\cdot) \right)\conjugate(y) \ge -j_0''(w(\omega);0)=0$ by definition of convex conjugation, the non-negativity of
$\parens*{\frac12 J''(w;\cdot)}\conjugate(\mu)$ follows.
\end{proof}

\begin{remark}
	\label{rem:conj}
As $j_0 \colon \R\to\R$, we have
\[
	j_0''(w(\omega);z) = j_0''(w(\omega);\sign(z)) z^2,
\]
which implies
\[
	\left(\frac12 j_0''(w(\omega);\cdot) \right)\conjugate(y) =
	\begin{cases}
	\frac12 \frac1{j_0''(w(\omega);\sign(y))}  y^2 & \text{ if }j_0''(w(\omega);\sign(y)) \ne 0,\\
	\delta_{\{0\}}(y) &\text{ otherwise.}
	\end{cases}
\]

\end{remark}

\subsection{Difference quotients of subdifferentials}
\label{subsec:dq_of_subdiff}
Next,
we study the limiting behaviour of a difference quotient of subdifferentials,
i.e.,
we show that
\begin{equation*}
	h_k \in \frac{ \partial J(w + t_k z) - \partial J(w)}{t_k}
\end{equation*}
converges in an appropriate way
towards
an element from
$\frac12 \partial J''(w, x; \cdot)(z)$,
where $x \in \partial J(w)$
under appropriate assumptions on $w$ and $z$.
Again,
we separate the smooth part in $J$ from the kinks.

\begin{lemma}
	\label{lem:dq_of_subdiff_max}
	Let $J$ be defined via \eqref{eq:integral_functionals_cont}
	with $j(w) := \max(w,0)$.
	Let $w\in C^1(\Omega)$ be given such that $\set{w=0}\cap \set{\nabla w=0}=\emptyset$
	and choose $x \in \partial J(w)$.
	Further, let $z \in C_c(\Omega)$ be arbitrary.
	Then, for every sequence $\seq{t_k} \subset \R^+$ with $t_k \searrow 0$,
	we can select
	a sequence $\seq{h_k} \subset \MM(\Omega)$
	with
	\begin{equation*}
		h_k \in
		\frac{ \partial J(w + t_k z) - x}{t_k}
		\qquad\forall k \in \N
	\end{equation*}
	such that
	\begin{equation}
		\label{eq:our_desire}
		h_k \weaklystar h,
		\qquad
		\norm{h_k}_{\MM(\Omega)} \to \norm{h}_{\MM(\Omega)},
	\end{equation}
	where
	$h \in \frac12\partial J''(w,x;\cdot)(z)$.
\end{lemma}
Note that the structural assumption
on $w$ implies that
$\partial J(w)$
and
$\partial J''(w,x;\cdot)(z)$
are singletons,
i.e., $x$ and $h$ are uniquely determined by $w$ and $z$.
Moreover, we even have $h_k \in L^1(\lambda)$,
but the limit is, in general, only a measure.

\begin{proof}
	We define
	\begin{equation*}
		x_k :=
		\chi_{\set{w + t_k z > 0}}
		+
		\chi_{\set{w + t_k z = 0, w < 0}}
	\end{equation*}
	and
	$h_k := (x_k - x)/t_k$.
	Then,
	it is clear that $x_k \in \partial J(w + t_k z)$
	and an easy distinction by cases shows
	\begin{equation*}
		h_k =
		\frac1{t_k} \sign(z) \chi_{\Omega_k}
		,
	\end{equation*}
	where $\Omega_k := \set{\sign(w) \ne \sign(w + t_k z)}$.
	Further, we note that the measure $h$ satisfies
	\begin{equation*}
		\dual{h}{\psi}
		=
		\int_{\set{w = 0}} \frac{\psi z}{\abs{\nabla w}} \d\HH^{d-1}
		\qquad
		\forall \psi \in C_0(\Omega),
	\end{equation*}
	see \eqref{eq:strong_second_order_subderivative}.
	In order to prove \eqref{eq:our_desire},
	we are going to apply \cref{thm:taylor_expansion_integral_rn}.
	For an arbitrary $\psi \in C_0(\Omega)$,
	we apply \eqref{eq:diff_int_rn_b}
	to obtain
	\begin{equation*}
		\dual{h_k}{\psi}
		=
		\frac1{t_k} \int_{\Omega_k} \psi \sign(z) \d\lambda
		\to
		\int_{\set{w = 0}} \frac{\psi z}{\abs{\nabla w}} \d\HH^{d-1}
		=
		\dual{h}{\psi}.
	\end{equation*}
	The convergence of norms follows from
	\begin{equation*}
		\norm{h}_{\MM(\Omega)}
		=
		\int_{\set{w = 0}} \frac{\abs{z}}{\abs{\nabla w}} \d\HH^{d-1},
		\quad
		\norm{h_k}_{\MM(\Omega)}
		=
		\norm{h_k}_{L^1(\lambda)}
		=
		\frac1{t_k} \int_{\Omega_k} 1 \d\lambda
	\end{equation*}
	and \eqref{eq:diff_int_rn_a}.
\end{proof}
Next, we address differentiable integrands.
To this end, we show
that the differentiability assumption on $j$ from \cref{lemma:zweeteableitung}
is equivalent to the directional differentiability of the first derivative.
\begin{lemma}
	\label{lem:quadratic_approximation_vs_second_derivative}
	Let $j \in C^1(\R)$ with locally Lipschitz continuous $j'$
	and $w,z \in \R$ be given.
	Then, the statements
	\begin{enumerate}
		\item
			\label{lem:quadratic_approximation_vs_second_derivative_1}
			$j''(w;z):=\lim_{t\searrow0} \frac2{t^2} (j(w+t z)- j(w) - t j'(w) z) \in \R$
			exists
		\item
			\label{lem:quadratic_approximation_vs_second_derivative_2}
			$D^2 j(w; z) := \lim_{t \searrow 0} \frac{j'(w + t z) - j'(w)}{t} \in \R$
			exists
	\end{enumerate}
	are equivalent
	and we have
	$j''(w; z) = D^2 j(w; z) z$.
\end{lemma}
\begin{proof}
	``\ref{lem:quadratic_approximation_vs_second_derivative_2}$\Rightarrow$\ref{lem:quadratic_approximation_vs_second_derivative_1}'':
	From \ref{lem:quadratic_approximation_vs_second_derivative_2}
	we have
	$j'(w + s z) - j'(w) = D^2 j(w; z) s + \oo(s)$ as $s \searrow 0$.
	Thus, the fundamental theorem of calculus implies
	\begin{equation*}
		j(w + t z) - j(w) - t j'(w) z
		=
		\int_0^t
		\parens[\big]{
			j'(w + s z)
			-
			j'(w)
		} z \d s
		=
		D^2 j(w; z) \frac{t^2}{2} z + \oo(t^2)
	\end{equation*}
	as $t \searrow 0$.
	This proves
	\ref{lem:quadratic_approximation_vs_second_derivative_1}
	with $j''(w; z) = D^2 j(w; z) z$.

	``\ref{lem:quadratic_approximation_vs_second_derivative_1}$\Rightarrow$\ref{lem:quadratic_approximation_vs_second_derivative_2}'':
	In case $z = 0$, the assertion is clear.
	It is sufficient to consider $z \ne 0$.
	We denote by $L$ the Lipschitz constant of $j'$ on $[w - \abs{z}, w + \abs{z}]$.
	Further, we define
	\[
	 q(t):=
	 j(w + t z)
	 -
	 j(w) - t j'(w) z - \frac{t^2}2 j''(w; z)
	 .
	\]
	Let $\varepsilon \in (0,1)$ be arbitrary.
	Due to \ref{lem:quadratic_approximation_vs_second_derivative_1},
	there exists $\hat t = \hat t_\varepsilon \in (0,1]$ such that
	\begin{equation*}
		\abs{ q(t) }
		\le
		\varepsilon t^2
		\qquad \forall t \in [0,\hat t].
	\end{equation*}
	For an arbitrary $t_1 \in (0,\hat t]$,
	we fix $t_0 := (1 - \sqrt{\varepsilon}) t_1 \in (0,t_1)$.
	Then we have
	\[
	 \frac{ q(t_1)-q(t_0)}{t_1-t_0}
	 =
	 \frac{j(w + t_1 z) - j(w + t_0 z)}{t_1 - t_0} - j'(w) z - t_1 j''(w; z) + \frac{t_1-t_0}2 j''(w; z)
	 .
	\]
	Using the Lipschitz continuity of $j'$, we get
	\begin{equation*}
		\abs*{ \frac{j(w+t_1z)-j(w+t_0z)}{t_1-t_0} - j'(w+t_1z)z}
		= \abs*{ \int_{t_0}^{t_1}\frac { j'(w+sz)z - j'(w+t_1z)z}{t_1-t_0}\d s}
		\le \frac L2 z^2 \parens{t_1-t_0}.
	\end{equation*}
Putting these expressions together, we find
\begin{align*}
	&
	\abs*{
		\bracks[\big]{j'(w + t_1 z) - j'(w)} z
		- t_1 j''(w; z)
	}
	\\&\qquad
	=
	\abs*{
		j'(w+t_1z)z - \frac{j(w+t_1z)-j(w+t_0z)}{t_1 - t_0}
		+
		\frac{q(t_1) - q(t_0)}{t_1 - t_0}
		- \frac{t_1-t_0}2 j''(w; z)
	}
	\\&\qquad
	\le  \frac L2 z^2 \parens{t_1-t_0}+  \varepsilon \frac{t_1^2 + t_0^2}{t_1 - t_0} + \frac{t_1 - t_0 }2 \abs{j''(w; z)}
	\le  \parens*{\frac L2 z^2 +  2 + \frac{1}2 \abs{j''(w; z)}}
	\sqrt{\varepsilon} t_1
	.
\end{align*}
Since $t_1 \in (0,\hat t_\varepsilon]$, $\epsilon \in (0,1)$ were arbitrary,
\ref{lem:quadratic_approximation_vs_second_derivative_2}
follows
with $D^2 j(w,z) = j''(w; z) / z$.
\end{proof}

\begin{remark}
	\label{rem:quadratic_approximation_vs_second_derivative}
	The statement
	\ref{lem:quadratic_approximation_vs_second_derivative_2}$\Rightarrow$\ref{lem:quadratic_approximation_vs_second_derivative_1}
	in \cref{lem:quadratic_approximation_vs_second_derivative}
	holds without the Lipschitz assumption on $j'$
	and this follows directly from the proof.
	The example
	$ j(s) = s^3 \sin(s^{-1}) $
	shows that the other direction may fail in absence of the Lipschitz assumption,
	since
	$ j'(s) = 3 s^2 \sin(s^{-1}) - s \cos(s^{-1}) $
	fails to be differentiable at $s = 0$.
\end{remark}

\begin{lemma}
	\label{lem:dq_of_subdiff_smooth}
	Let $J$ be defined via \eqref{eq:integral_functionals_cont}
	with $j \in C^1(\R)$ with locally Lipschitz continuous $j'$,
	such that one of the assumptions of \cref{lem:quadratic_approximation_vs_second_derivative}
	is satisfied.
	Let $w\in C^1(\Omega)$ be bounded
	and set $x = J'(w)$.
	Further, let $z \in C_0(\Omega)$ be arbitrary.
	Then, for every sequence $\seq{t_k} \subset \R^+$ with $t_k \searrow 0$,
	the sequence $\seq{h_k} \subset \MM(\Omega)$
	with
	\begin{equation*}
		h_k :=
		\frac{ J'(w + t_k z) - x}{t_k}
		\qquad\forall k \in \N
	\end{equation*}
	satisfies
	\begin{equation}
		\label{eq:our_desire_again}
		h_k \weaklystar h,
		\qquad
		\norm{h_k}_{\MM(\Omega)} \to \norm{h}_{\MM(\Omega)},
	\end{equation}
	where $h \in \MM(\Omega)$ is defined via
	\begin{equation*}
		\dual{h}{\psi}
		=
		\int_\Omega D^2 j(w; z) \psi \d\lambda
		\qquad\forall \psi \in C_0(\Omega)
		.
	\end{equation*}
\end{lemma}
\begin{proof}
	We take $x_k = J'(w + t_k z)$ and have
	\begin{equation*}
		\dual{x_k}{\psi} = \int_\Omega j'(w + t_k z) \psi \d\lambda,
		\qquad
		\dual{x  }{\psi} = \int_\Omega j'(w        ) \psi \d\lambda.
	\end{equation*}
	Thus,
	\begin{equation*}
		\dual{h_k}{\psi} =
		\int_\Omega \frac{j'(w + t_k z) - j'(w)}{t_k} \psi \d\lambda,
		\qquad
		\norm{h_k}_{\MM(\Omega)}
		=
		\int_\Omega \abs*{\frac{j'(w + t_k z) - j'(w)}{t_k}} \d\lambda.
	\end{equation*}
	By dominated convergence,
	we get $h_k \weaklystar h$ and $\norm{h_k}_{\MM(\Omega)} \to \norm{h}_{\MM(\Omega)}$.
\end{proof}
In the setting of \cref{lem:dq_of_subdiff_smooth},
we even have $h, h_k \in L^1(\lambda)$.

Now, we combine \cref{lem:dq_of_subdiff_max,lem:dq_of_subdiff_smooth}.
\begin{theorem}
	\label{thm:dq_of_subdiff}
	We consider the setting of \cref{asm:dzweij_pw}
	with a convex $j_0$.
	Let a bounded $w \in C^1(\Omega)$ be given such that \eqref{eq_ass_w_pw} holds
	and choose $x \in \partial J(w)$.
	Further, let $z \in C_c(\Omega)$ be arbitrary.
	Then, for every sequence $\seq{t_k} \subset \R^+$ with $t_k \searrow 0$,
	we can select
	a sequence $\seq{h_k} \subset \MM(\Omega)$
	with
	\begin{equation}
		\label{eq:FD_inclusion}
		h_k \in
		\frac{ \partial J(w + t_k z) - x}{t_k}
		\qquad\forall k \in \N
	\end{equation}
	such that
	\begin{equation}
		\label{eq:our_desire_yet_another}
		h_k \weaklystar h,
		\qquad
		\norm{h_k}_{\MM(\Omega)} \to \norm{h}_{\MM(\Omega)},
	\end{equation}
	where
	$h \in \frac12\partial J''(w,x;\cdot)(z)$,
	i.e.,
	\begin{equation}
		\label{eq:subdif_j_pp}
		\dual{h}{\psi}
		=
		\int_\Omega \frac12 \partial j_0''(w; \cdot)(z) \psi \d\lambda
		+
		\sum_{i=1}^m
		2a_i \int_{\set{w=b_i}} \frac{ z \psi }{\abs{\nabla w}} \d\HH^{d-1}
		\qquad\forall \psi \in C_0(\Omega).
	\end{equation}
\end{theorem}
Note that the second integral in \eqref{eq:subdif_j_pp}
is well defined since $w$ satisfies \eqref{eq_ass_w_pw} and since $z$ has a compact support.
\begin{proof}
	We define the functionals
	\begin{equation*}
		J_0(\tilde w) := \int_\Omega j(\tilde w) - \sum_{i=1}^m a_i(\tilde w-b_i) \d\lambda,
		\qquad
		J_i(\tilde w) := \int_\Omega \max(\tilde w - b_i, 0) \d\lambda
		\qquad\forall \tilde w \in C_0(\Omega),
	\end{equation*}
	which implies with $a_0:=\frac12$
	\[
	 J(\tilde w) = \sum_{i=0}^m 2a_i J_i(\tilde w).
	\]
	Due to \eqref{eq_ass_w_pw},
	these functionals are differentiable at $w$
	and we select the unique elements $x^i \in \partial J_i(w)$.
	From \cref{lem:dq_of_subdiff_max,lem:dq_of_subdiff_smooth}
	we obtain sequences $\seq{h_k^i}$
	with
	$h_k^i \in \frac{\partial J_i(w + t_k z) - x^i}{t_k}$
	and
	$h_k^i \weaklystar h^i$, $\norm{h_k^i}_{\MM(\Omega)} \to \norm{h^i}_{\MM(\Omega)}$,
	where
	\begin{align*}
		\dual{h^0}{\psi}
		&=
		\int_\Omega D^2 j_0(w; z) \psi \d\lambda
		&&\forall \psi \in C_0(\Omega),
		\\
		\dual{h^i}{\psi}
		&=
		\int_{\set{w=b_i}} \frac{ z \psi }{\abs{\nabla w}} \d\HH^{d-1}
		&&\forall \psi \in C_0(\Omega), i = 1,\ldots, m.
	\end{align*}
	From the relation $j_0''(w; z) = D^2 j_0(w; z) z$,
	we get $\frac12 \partial j_0''(w; \cdot)(z) = D^2 j_0(w; z)$.
	This implies
	\begin{equation*}
		h^0
		=
		\frac12 \partial J_0''(w; \cdot)(z),
	\end{equation*}
	see \cref{thm:rockit}.
	The sum rule of convex analysis implies that $h_k := \sum_{i = 0}^m 2 a_i h_k^i$
	satisfies \eqref{eq:FD_inclusion}.
	The first limit relation in \eqref{eq:our_desire_yet_another} is clear.
	The second relation follows from
	\begin{align*}
		\norm{h}
		&
		\le \liminf_{k \to \infty} \norm{h_k}
		\le \limsup_{k \to \infty} \norm{h_k}
		\le \sum_{i = 0}^m \lim_{k \to \infty} 2 a_i \norm{h_k^i}
		=
		\sum_{i = 0}^m 2 a_i \norm{h^i}
		=
		\norm{h}.
	\end{align*}
	In the last equality, we have used the precise structure of $h^i$.
\end{proof}
For a similar result in Hilbert spaces,
we refer to
\cite[Theorem~3.9]{Do}.
Therein,
it is shown that the strong twice epi-differentiability of a functional
$J$
is equivalent to
$\partial J$ being proto-differentiable with maximally monotone proto-derivative.
The proof utilizes Attouch's theorem on the Mosco-convergence of functionals
and this is not available in non-reflexive spaces,
see \cite[Theorem~2.5]{CombariThibault1998} and the comment thereafter.

\subsection{A localized descent lemma}
\label{subsec:descent}
As a last preparation,
we state a localized descent lemma for $J$
at a point $w$ satisfying a structural assumption.

\begin{lemma}
	\label{lem:descent}
 We consider the setting of \cref{asm:dzweij_pw}.
	Further, let $w \in C^1(\Omega) \cap C_0(\Omega)$ be given such that
	there exist $C, \eta > 0$ with
	\begin{equation}
		\label{eq:structural}
		\sum_{i = 1}^m \lambda\parens{\set{ \abs{ w - b_i } \le \varepsilon }}
		\le
		C \varepsilon
		\qquad
		\forall \varepsilon \in [0,\eta].
	\end{equation}
	Then, $J$ is (Gâteaux) differentiable at $w$ and we set $x = J'(w)$.
	Moreover, there exists a constant $\Lambda > 0$
	such that
	\begin{equation}
		\label{eq:descent}
		J(v)
		\le
		J(w) + \dual{x}{v - w}
		+
		\frac \Lambda 2 \norm{v - w}_{C_0(\Omega)}^2
		\qquad\forall v \in C_0(\Omega):\ \norm{v - w}_{C_0(\Omega)} \le \eta
		.
	\end{equation}
\end{lemma}
\begin{proof}
	Due to \eqref{eq:structural}, we have
	\begin{equation*}
		\sum_{i = 1}^m \lambda\parens{\set{ w = b_i }}
		=
		0
	\end{equation*}
	and this is enough to ensure that $\partial J(w)$ is a singleton.
	Consequently, $J$ is (Gâteaux) differentiable at $w$.

	Next, it is easy to see that it is sufficient to validate \eqref{eq:descent}
	for the components
	\begin{equation*}
		J_0(v) := \int_\Omega j_0(v(\omega)) \d\lambda(\omega),
		\qquad
		J_i(v) := \int_\Omega \abs{v - b_i} \d\lambda.
	\end{equation*}
	Since $j_0'$ is assumed to be locally Lipschitz,
	there exists a constant $L \ge 0$
	such that
	\begin{equation*}
		\abs{ j_0'(r) - j_0'(s) }
		\le
		L \abs{r - s}
		\qquad
		\forall r,s \in [-T,T]
	\end{equation*}
	with $T := \norm{w}_{C_0(\Omega)} + \eta$.
	Let $v \in C_0(\Omega)$ with $\norm{v - w}_{C_0(\Omega)} \le \eta$
	be given.
	The Lipschitz estimate implies
	\begin{equation*}
		\abs{ \dual{J_0'(v) - J_0'(w)}{z} }
		\le
		\int_\Omega \abs{j_0'(v(\omega)) - j_0'(w(\omega))} \abs{z} \d\lambda(\omega)
		\le
		L \lambda(\Omega) \norm{ v - w }_{C_0(\Omega)} \norm{ z }_{C_0(\Omega)}.
	\end{equation*}
	From \cref{thm:rockit}
	we know that $J_0$ is (Gâteaux) differentiable on $C_0(\Omega)$.
	Hence, the fundamental theorem of calculus implies
	\begin{align*}
		J_0(v) - J_0(w)
		-\dual{J_0'(w)}{v - w}
		&=
		\int_0^1 \dual{J_0'( w + t \, (v-w) ) - J_0'(w)}{v - w}_{C_0(\Omega)} \d t
		\\
		&\le
		L \lambda(\Omega) \int_0^1 t \norm{v - w}_{C_0(\Omega)}^2 \d t
		=
		\frac {L \lambda(\Omega)} 2 \norm{w - v}_{C_0(\Omega)}^2.
	\end{align*}
	Next, we consider $J_i$ for $i = 1,\ldots, m$.
	By using the identities
	\begin{align*}
		\abs{r}
		&=
		\abs{s}
		+
		\sign(s) (r - s)
		&&\text{if } \sign(r) = \sign(s)
		\\
		\abs{ r }
		&\le
		\abs{ s }
		+
		\sign(s) (r - s)
		+
		2 \abs{r - s}
		&&\text{if } \sign(r) \ne \sign(s)
	\end{align*}
	for all $r,s \in \R$ with $s \ne 0$,
	we obtain
	\begin{align*}
		J_i(v)
		-
		J_i(w) - \dual{J_i'(w)}{v - w}
		&\le
		\int_{\set{\sign(v-b_i) \ne \sign(w-b_i)}} 2 \abs{ v - w } \d\lambda
		\\&\le
		2 \lambda(\set{\sign(v-b_i) \ne \sign(w-b_i)}) \norm{ v - w }_{C_0(\Omega)}.
	\end{align*}
	On $\set{\sign(v-b_i) \ne \sign(w-b_i)}$ we have $\abs{v - w} \ge \abs{w-b_i}$.
	Thus,
	\begin{align*}
		\lambda(\set{\sign(v-b_i) \ne \sign(w-b_i)})
		&\le
		\lambda( \set{\abs{w-b_i} \le \abs{v - w}})
		\\&\le
		\lambda( \set{\abs{w-b_i} \le \norm{v - w}_{C_0(\Omega)}})
		\le
		C \norm{v - w}_{C_0(\Omega)}
	\end{align*}
	for $\norm{v - w}_{C_0(\Omega)} \le \eta$.
	This shows
	\begin{align*}
		J_i(v)
		\le
		J_i(w) + \dual{J_i'(w)}{v - w}
		+
		2 C \norm{ v - w }_{C_0(\Omega)}^2
		.
	\end{align*}
	Combining the estimates for $J_0$ and $J_i$ yields the claim.
\end{proof}

A condition similar to \eqref{eq:structural}
was first used (almost simultaneously and independently)
in
\cite{DeckelnickHinze2012,WachsmuthWachsmuth2009}.
This condition is related to \eqref{eq_int_nablaw_z_finite}, i.e., $\frac1{\abs{\nabla w}}|_\ZZ \in L^1(\HH^{d-1}|_\ZZ)$,
see \cref{lem_measure_condition} below.

\begin{remark}
	\label{rem:on_structural}
	The condition \eqref{eq:structural}
	follows almost from \eqref{eq_ass_w_pw}.
	Indeed, let us assume that $w \in C^1(\Omega)$
	and its derivatives can be extended continuously to $\bar\Omega$,
	i.e., $w \in C^1(\bar\Omega)$
	and that
	\begin{equation*}
		w(\omega) = b_i
		\quad\Rightarrow\quad
		\nabla w(\omega) \ne 0
		\qquad
		\forall w \in \bar\Omega.
	\end{equation*}
	Note that \eqref{eq_ass_w_pw} is the same condition, but only for all $\omega \in \Omega$.
	Then, \cite[Lemma~3.2]{DeckelnickHinze2012}
	shows that \eqref{eq:structural} is satisfied.
\end{remark}
\begin{remark}
	\label{rem:global_lipschitz}
	Condition \eqref{eq:structural} implies
	\begin{equation*}
		\sum_{i = 1}^m \lambda\parens{\set{ \abs{ w - b_i } \le \varepsilon }}
		\le
		\max\set{C, m \lambda(\Omega) / \eta} \varepsilon
		\qquad
		\forall \varepsilon \ge 0.
	\end{equation*}
	If we assume improved regularity of $j_0$, i.e., if $j_0'$ is globally Lipschitz,
	we can show that \eqref{eq:descent} holds even for $\eta = \infty$.
\end{remark}

Let us discuss how condition \eqref{eq:structural} is related to \eqref{eq_int_nablaw_z_finite}, i.e., $\frac1{\abs{\nabla w}}|_\ZZ \in L^1(\HH^{d-1}|_\ZZ)$.

\begin{lemma}\label{lem_measure_condition}
Let $w \in C^1(\Omega)$ satisfy \eqref{eq_ass_w_pw}. Then \eqref{eq:structural} implies \eqref{eq_int_nablaw_z_finite}.
In fact,
\[
	2 \int_\ZZ \frac1{\abs{\nabla w}} \d\HH^{d-1} \le C
\]
with $C$ from  \eqref{eq:structural}
and $\ZZ$ as in \eqref{eq:def_ZZ}.
In addition, \eqref{eq_int_nablaw_z_finite} implies
\[
	\limsup_{t\searrow0} \frac1t \sum_{i=1}^m \lambda\parens{K \cap \set{ \abs{ w - b_i } \le t}}
	\le
	2 \int_\ZZ \frac1{\abs{\nabla w}} \d\HH^{d-1},
\]
for all compact $K\subset \Omega$.
\end{lemma}
\begin{proof}
Assume \eqref{eq:structural}.
	Let $t\in (0,\eta]$. Let $z\in C_c(\Omega)$ with  $\norm{z}_{C_0(\Omega)}\le 1$ be given.
	Define
	\[\begin{aligned}
        \Omega_t^+ & := \{ \sign(w-b_i) \ne \sign(w-b_i +tz ) \},\\
        \Omega_t^- & := \{ \sign(w-b_i) \ne \sign(w-b_i -tz ) \}.
        \end{aligned}
	\]
	Then we find
	\[
        \Omega_t := \Omega_t^+ \cup \Omega_t^- \subset \{ \abs{ w - b_i } \le t \}.
	\]
	By \eqref{eq:structural}, we get $\lambda(\Omega_t) \le C t$.
	Then \cref{thm:taylor_expansion_integral_rn} \eqref{eq:diff_int_rn_a} implies
	\[
		C  \ge \frac1t \lambda(\Omega_t^+ \cup \Omega_t^-)= \frac1t\int_{\Omega_t^+ \cup \Omega_t^-} \dx \to 2\int_{\{w=b_i\}} \frac{\abs{z}} {\abs{\nabla w}} \d\HH^{d-1}.
	\]
	This implies $2 \int_{\{w=b_i\}} \frac1 {\abs{\nabla w}} \d\HH^{d-1} \le C$, and \eqref{eq_int_nablaw_z_finite} follows.

Let \eqref{eq_int_nablaw_z_finite} be satisfied. Take $K\subset \Omega$ compact.
Let $z\in C_c(\Omega)$ with $\norm{z}_{C_0(\Omega)}\le 1$ such that $z=1$ on $K$.
Then one easily verifies
	\[
	        \{ \abs{ w - b_i } \le t \} \cap K \subset \Omega_t^+ \cup \Omega_t^- .
	\]
Together with \cref{thm:taylor_expansion_integral_rn} \eqref{eq:diff_int_rn_a}  we obtain
\begin{align*}
	\limsup_{t\searrow0} \frac1t \lambda\parens{K \cap \set{ \abs{ w - b_i } \le t}}
	&\le
	\lim_{t\searrow0} \frac1t \int_{\Omega_t^+ \cup \Omega_t^-} \d\lambda
	\
	=
	2 \int_{\set{w=b_i}}\frac{\abs{z}}{\abs{\nabla w}} \d\HH^{d-1}
	\\&\le
	2 \int_{\set{w=b_i}} \frac1{\abs{\nabla w}} \d\HH^{d-1},
\end{align*}
which proves the claim.
\end{proof}

\begin{corollary}
Let $w \in C^1(\Omega) \cap C_0(\Omega)$ satisfy \eqref{eq_ass_w_pw}.
Assume $b_i \ne 0$ for all $i$.
Then
	\eqref{eq:structural} and \eqref{eq_int_nablaw_z_finite} are equivalent.
	In this case
	\[
	\limsup_{t\searrow0} \frac1t \sum_{i=1}^m \lambda\parens{ \set{ \abs{ w - b_i } \le t}}
	=
	2 \int_\ZZ \frac1{\abs{\nabla w}} \d\HH^{d-1}
		.
	\]
\end{corollary}
\begin{proof}
Assume  \eqref{eq_int_nablaw_z_finite}.
	The assumptions on $b_i$ imply that $\ZZ$ (as in \eqref{eq:def_ZZ}) has a positive distance to the boundary $\partial \Omega$,
	and $\bigcup_{i=1}^m \set{ \abs{ w - b_i } \le t}$ is a compact subset of $\Omega$ for $t>0$ small enough.
Then \cref{lem_measure_condition} implies
\[
	\limsup_{t\searrow0} \frac1t \sum_{i=1}^m \lambda\parens{ \set{ \abs{ w - b_i } \le t}}
	\le
	2 \int_\ZZ \frac1{\abs{\nabla w}} \d\HH^{d-1},
\]
which in turn implies \eqref{eq:structural}.
\end{proof}

\section{Second order derivatives of integral functionals over measures}
\label{sec:sod_integral_functionals_measures}
In this section,
we are going to study
integral functionals over the space of measures
$\MM(\Omega)$.
Recall that $\MM(\Omega)$
is the (topological) dual space of $C_0(\Omega)$.
In contrast, $L^1(\lambda)$
cannot be a dual space
(unless $\Omega = \emptyset$)
and this is the reason for working in $\MM(\Omega)$.

The space $L^1(\lambda)$ is considered as a subspace of $\MM(\Omega)$
in the sense that
a function $u \in L^1(\lambda)$ is identified with the measure
$u \lambda \in \MM(\Omega)$ defined via
\begin{equation*}
	(u \lambda)(A)
	\mapsto
	\int_A u \d\lambda
	\qquad\forall A \in \BB(\Omega).
\end{equation*}
Now, let the integral functional $G$ be induced by a convex and lower semicontinuous function $g \colon \R \to \bar\R$, i.e.,
\begin{equation*}
	G(u) := \int_\Omega g(u(\omega)) \d\lambda(\omega)
	\qquad\forall u \in L^1(\lambda).
\end{equation*}
Note that $G(u) \in (-\infty,\infty]$ is well defined
since $g$ possesses an affine minorant and $u \in L^1(\lambda)$.
We extend this functional to $\MM(\Omega)$ via
\begin{equation*}
	G(\mu)
	:=
	\begin{cases}
	\int_\Omega g\parens[\big]{ \frac{\d\mu}{\d\lambda}(\omega) } \d\lambda(\omega)
	&\text{if } \mu \ll \lambda, \\
	+\infty & \text{otherwise}
	\end{cases}
	\qquad
	\forall \mu \in \MM(\Omega).
\end{equation*}
This definition coincides with the original definition of $G$
on $L^1(\lambda)$
since
$u\lambda \ll \lambda$ and
\begin{equation*}
	\frac{\d (u \lambda)}{\d\lambda} = u
	\qquad\forall u \in L^1(\lambda).
\end{equation*}
Further,
$\dom(G) \subset L^1(\lambda)$, since
all $\mu \in \MM(\Omega)$
are finite measures, and hence $\d\mu/\!\d\lambda \in L^1(\lambda)$.
Thus, we have extended $G$ by $\infty$
on the set $\MM(\Omega) \setminus L^1(\lambda)$.

The aim of this section
is the computation of the second subderivative
of $G$ under some structural assumptions
and the verification of its strict twice epi-differentiability.
This will be done by deriving
estimates from below and from above.

In order to apply the results from \cref{sec:sod_integral_functionals_cont},
we restrict ourselves to functions $g$
whose conjugates fit into the setting of the previous section.
\begin{assumption}
	\label{asm:structure_of_conjugate}
	We assume that the conjugate $g\conjugate =: j$
	satisfies the assumptions posed on $j$
	in \cref{asm:dzweij_pw}.
	Moreover, for each bounded set $W\subset \R$ there exists a constant $C_j \ge 1$ such that the smooth part $j_0$
	satisfies
	\begin{equation*}
		j_0''(w; -1) = 0
		\text{ or }
		j_0''(w; +1) = 0
		\text{ or }
		C_j^{-1} j_0''(w; -1)
		\le
		j_0''(w; +1)
		\le
		C_j j_0''(w; -1)
	\end{equation*}
	for all $w \in W$.
\end{assumption}
We require that \cref{asm:structure_of_conjugate}
holds throughout this \cref{sec:sod_integral_functionals_measures}.

\begin{remark}
	\label{rem:not_restrictive}
	The second part of \cref{asm:structure_of_conjugate} is not very restrictive
	and is satisfied by many reasonable functions $j_0$.
	However, it is possible to construct some artificial counterexamples:
	We choose $\zeta \in L^\infty(\R)$ with
	\begin{equation*}
		\zeta(t) =
		\begin{cases}
			1/n & \text{if } t \in (1/n, 1/n+2^{-n}) \text{ for some } n \in \N, \\
			1 & \text{else}.
		\end{cases}
	\end{equation*}
	We set
	\begin{equation*}
		j(t) = \int_0^t \int_0^s \zeta(r) \d r \d s.
	\end{equation*}
	It is clear that $j$ is convex and differentiable, and that $j'$ is globally Lipschitz.
	With \cref{lem:quadratic_approximation_vs_second_derivative}
	one can check that
	\begin{align*}
		j''(t; 1)
		&=
		1/n \quad\text{if } t \in [1/n, 1/n+2^{-n}) \text{ for some } n \in \N,
		\\
		j''(t; -1)
		&=
		1/n \quad\text{if } t \in (1/n, 1/n+2^{-n}] \text{ for some } n \in \N,
	\end{align*}
	and in all other cases we have $j''(t; \pm 1) = 1$.
	Thus, the second part of \cref{asm:structure_of_conjugate} is violated.
\end{remark}

As a preliminary result,
we give a functional
$J \colon C_0(\Omega) \to \R$
such that the conjugate of $J$ is $G$.
The computation of conjugate functions of integral functionals on continuous functions
is addressed in \cite[Section~3]{Rockafellar1971}, see in particular \cref{thm:rockit}.
Therefore,
we expect that $J$ can be obtained by integrating over the convex conjugate of $g$.
\begin{lemma}
	\label{lem:preconjugate}
	We define
	$J \colon C_0(\Omega) \to \R$ via
	\begin{equation*}
		J(z)
		:=
		\int_\Omega g\conjugate(z(\omega)) \d\lambda(\omega).
	\end{equation*}
	Then, the conjugate of $J$ is given by $G$.
\end{lemma}
\begin{proof}
	This directly follows from \cref{thm:rockit}
	with the choices
	$\iota(\omega, x) = g\conjugate(x)$
	and
	$\nu = \lambda$.
\end{proof}

It would be tempting to obtain the weak-$\star$ second subderivative
by computing the second subderivative of $J$
and a convex conjugation argument.
For strongly twice epi-differentiable functions
defined on reflexive Banach spaces,
this is indeed possible, see \cite[Theorem~2.5]{Do}.
The argument in \cite[Theorem~2.5]{Do}
uses the Wijsman--Mosco theorem on the continuity of the
Legendre--Fenchel transform w.r.t.\ Mosco epi-convergence of convex functions,
i.e.\ the Mosco epi-convergence of $\varphi_n$ to $\varphi$
is equivalent to the Mosco epi-convergence $\varphi_n\conjugate$ to $\varphi\conjugate$
in reflexive spaces.
In non-reflexive spaces, this theorem fails,
see the discussion in \cite[Section~1]{Zabell1992}
and \cite[Theorem~3.3]{BeerBorwein1990}.
Hence, we have to use a direct argument
which connects the weak-$\star$ second subderivative of $G$
with the strong second subderivative of $J$.

\subsection{Estimates from below}

In this section, we are going to give a lower estimate for the second subderivative of $G$.
The arguments are similar to those of \cite[Lemma~6.8, Proposition~6.9]{ChristofWachsmuth2017:1}.
To this end, we utilize the preconjugate function $J$.
Note that the next lemma also works in a more general setting.
Therefore, we state it with $Y$ and $X$ instead of $C_0(\Omega)$ and $\MM(\Omega)$.
\begin{lemma}
	\label{lem:lower_bound}
	Let $J \colon Y\to \bar \R$ be convex, proper, lower semicontinuous  such that $J\conjugate= G$.
	Let $x \in \dom(G)$, $w \in \partial G(x) \cap Y$ be fixed
	such that
	$J$ is strongly-strongly twice epi-differentiable at $w$ for $x$.
	Then,
	\begin{equation}
		\label{eq:lower_estimate_2}
		\frac12 G''(x, w; h)
		\ge
		\parens*{\frac12 J''(w,x;\cdot)}\conjugate(h)
		\qquad\forall h \in X
		.
	\end{equation}
\end{lemma}
\begin{proof}
	Let $z \in Y$ be arbitrary.
Further, let $\seq{t_k} \subset (0,\infty)$ and $\seq{h_k} \subset X$
be arbitrary sequences with
$t_k \searrow 0$ and $h_k \weaklystar h$ in $X$.
Next, we choose a sequence $\seq{z_k}$ in $Y$  with $z_k \to z$
depending on $\seq{t_k}$, see \eqref{eq:conv_second_order_J} below.
Since $G$ is the convex conjugate of $J$, we have
the Fenchel--Moreau inequality
\begin{equation}
	\label{eq:FMI}
	G(x + t_k h_k) \ge \dual{x + t_k h_k}{w + t_k z_k} - J(w + t_k z_k)
\end{equation}
and,
due to $w \in \partial G(x)$,
we find
\begin{equation*}
	G(x) = \dual{x}{w} - J(w)
	.
\end{equation*}
Taking the difference of these two relations and dividing by $t_k^2/2$ yields
\begin{equation*}
	\frac{G(x + t_k h_k) - G(x) - t_k \dual{h_k}{w}}{t_k^2/2}
	\ge
	-\frac{J(w + t_k z_k) - J(w) - t_k \dual{x}{z_k}}{t_k^2/2}
	+
	2 \dual{h_k}{z_k}.
\end{equation*}
Due to the strongly-strongly twice epi-differentiability of $J$,
the sequence $\seq{z_k}$ can be chosen such that
we have the convergence
\begin{equation}
	\label{eq:conv_second_order_J}
	J''(w,x;z)
	=
	\lim_{k \to \infty} \frac{J(w + t_k z_k) - J(w) - t_k \dual{x}{z_k}}{t_k^2/2}
	\in [0,\infty]
	.
\end{equation}
Then, we can pass to the limit in the above inequality and obtain
\begin{equation*}
	\liminf_{k \to \infty}\frac{G(x + t_k h_k) - G(x) - t_k \dual{h_k}{w}}{t_k^2/2}
	\ge
	-J''(w,x;z)
	+
	2 \dual{h}{z}
	.
\end{equation*}
Since the right-hand side is independent of the sequences $\seq{h_k}$ and $\seq{t_k}$,
we can take the infimum over all sequences $t_k \searrow 0$ and $h_k \weaklystar h$ in $X$.
Thus,
\begin{equation}
	\label{eq:lower_estimate}
	\frac12 G''(x, w; h)
	\ge
	-\frac12 J''(w,x;z)
	+
	\dual{h}{z}.
\end{equation}
Since $z \in Y$ was arbitrary,
we can take the supremum w.r.t.\ $z \in Y$
in the right-hand side.
This yields
\eqref{eq:lower_estimate_2}.
\end{proof}
We mention that
the strongly-strongly twice epi-differentiability of $J$ was investigated in \cref{thm_dzweij_pw},
and
an expression for the conjugate
of the strong second subderivative of $J$
is given in \cref{thm:conj_subderivative}.

\subsection{Strictly twice epi-differentiability of \texorpdfstring{$G$}{G}}
\label{subsec:upper_bound}

Next, we give an upper bound
by considering sequences $\seq{z_k}$ and $\seq{h_k}$
for which we have equality in \eqref{eq:FMI}.
This is possible due to \cref{thm:dq_of_subdiff}.

\begin{lemma}
	\label{lem:upper_bound}
	Let $x \in \dom(G)$, $w \in \partial G(x) \cap C_0(\Omega)$ be fixed
	such that
	$w$ belongs to
	$C^1(\Omega)$
	and satisfies
	\eqref{eq_ass_w_pw}.
	For an arbitrary
	$z \in C_c(\Omega)$,
	we select
	\begin{equation*}
		h \in \frac12 \partial J''(w, x; \cdot) (z)
		.
	\end{equation*}
	Then,
	\begin{equation*}
		\frac12 G''(x, w; h)
		=
		\parens*{\frac12 J''(w,x;\cdot)}\conjugate(h)
		.
	\end{equation*}
	Moreover,
	for every sequence $\seq{t_k}$ with $t_k \searrow 0$,
	we can select a sequence $\seq{h_k} \subset \MM(\Omega)$
	such that $h_k \weaklystar h$, $\norm{h_k}_{\MM(\Omega)} \to \norm{h}_{\MM(\Omega)}$
	and
	\begin{equation*}
		G''(x, w; h)
		=
		\lim_{k \to \infty}
		\frac{G(x + t_k h_k) - G(x) - t_k \dual{h_k}{w}}{t_k^2/2}
		.
	\end{equation*}
\end{lemma}
We mention that $h$ is uniquely determined by $z$,
cf.\ \eqref{eq:subdif_j_pp}.
\begin{proof}
For an arbitrary sequence $t_k \searrow 0$,
we
use \cref{thm:dq_of_subdiff} to select
\begin{equation*}
	h_k \in \frac{\partial J(w + t_k z) - x}{t_k}.
\end{equation*}
with
$h_k \weaklystar h$, $\norm{h_k}_{\MM(\Omega)} \to \norm{h}_{\MM(\Omega)}$.
Due to $x + t_k h_k \in \partial J(w + t_k z)$,
we have
\begin{align*}
	J(w + t_k z) + G(x + t_k h_k) &= \dual{x + t_k h_k}{w + t_k z},
	\\
	\mrep{J(w)}{J(w + t_k z)} + \mrep{G(x)}{G(x + t_k h_k)} &= \dual{x}{w}.
\end{align*}
Taking the difference and dividing by $t_k^2$ yields
\begin{equation*}
	\frac12 \frac{G(x + t_k h_k) - G(x) - t_k \dual{h_k}{w}}{t_k^2/2}
	+
	\frac12 \frac{J(w + t_k z) - J(w) - t_k \dual{x}{z}}{t_k^2/2}
	=
	\dual{h_k}{z}.
\end{equation*}
The right-hand side converges towards $\dual{h}{z}$.
Via \cref{lemma:constantin,lemma:zweeteableitung},
the second addend on the left-hand side converges towards $J''(w;z)$.
Thus, we have
\begin{align*}
	\frac12 G''(x, w; h) + \frac12 J''(w; z)
	&\le
	\frac12 \liminf_{k \to \infty}
	\frac{G(x + t_k h_k) - G(x) - t_k \dual{h_k}{w}}{t_k^2/2}
	+
	\frac12 J''(w; z)
	\\
	&\le
	\frac12 \limsup_{k \to \infty}
	\frac{G(x + t_k h_k) - G(x) - t_k \dual{h_k}{w}}{t_k^2/2}
	+
	\frac12 J''(w; z)
	\\&\le
	\dual{h}{z}
	.
\end{align*}
Now, \eqref{eq:lower_estimate}
implies that all three inequalities are satisfied with equality
and together with \eqref{eq:lower_estimate_2} this shows the claim.
\end{proof}
We remark that the arguments in the proof of \cref{lem:upper_bound}
can also be utilized in a more general situation
if the preconjugate $J$
satisfies the assertion of \cref{lemma:constantin,lemma:zweeteableitung,thm:dq_of_subdiff}
with $C_c(\Omega)$ replaced by some subspace of $Y$.

In order to prove the equality
\begin{equation*}
    \frac12 G''(x, w; h)
    =
    \parens*{\frac12 J''(w,x;\cdot)}\conjugate(h)
		\qquad\forall h \in \MM(\Omega),
\end{equation*}
we will utilize the density argument  \cite[Lemma~3.2]{ChristofWachsmuth2019}.
Therefore, we need to show that for each $h\in \MM(\Omega)$ there exists a sequence $(h_k)$ with $h_k \in \frac12 \partial J''(w,x; \cdot)(z_k)$
for some $z_k \in C_c(\Omega)$ such that
\[
 \parens*{\frac12 J''(w,x;\cdot)}\conjugate(h_k) \to
 \parens*{\frac12 J''(w,x;\cdot)}\conjugate(h) .
\]
In order to construct this sequence, we will use the following approximation result.

\begin{lemma}
	\label{lem:cts}
Let $\nu$ be a positive Borel measure on $\Omega$ such that $\nu(K)<\infty$ for all compact subsets $K \subset \Omega$.
Let $u\in L^1(\nu) \cap L^2(\nu)$ be given.
Then there exists a sequence $\seq{u_k} \subset C_c(\Omega)$ such that
$u_k \to u$ in $L^1(\nu) \cap L^2(\nu)$ and pointwise almost everywhere.
In addition, there exists $g\in L^1(\nu) \cap L^2(\nu)$ such that $\abs{u_k} \le g$.
\end{lemma}
\begin{proof}
By \cite[Theorem 2.18]{Rudin1987}, $\nu$ is regular.
Suppose $u\ge0$.
Due to \cite[Theorem 3.14]{Rudin1987} there are sequences $\seq{v_k}, \seq{w_k} \subset C_c(\Omega)$ such that
$v_k \to u$ in $L^1(\nu)$ and $w_k \to u$ in $L^2(\nu)$.
W.l.o.g.\ we can assume $v_k \ge0$, $w_k\ge0$, $v_k \to u$ and $w_k \to u$ pointwise almost everywhere and
the existence of $g_1\in L^1(\nu)$, $g_2\in L^2(\nu)$ with $v_k\le g_1$, $w_k\le g_2$.
Define $u_k := \min(v_k,w_k)$. Then $u_k \le \min(g_1,g_2)=:g \in L^1(\nu) \cap L^2(\nu)$
and $u_k \to u$ in $L^1(\nu) \cap L^2(\nu)$ follows by the dominated convergence theorem.
For general $u$, we can apply this construction to the positive and negative part of $u$.
\end{proof}

\begin{lemma}
\label{thm:recovery_seq}
	Let $x \in \dom(G)$, $w \in \partial G(x) \cap C_0(\Omega)$ be fixed
	such that
	$w$ belongs to
	$C^1(\Omega)$
	and \eqref{eq_int_nablaw_z_finite} is satisfied.
Let $h\in \MM(\Omega)$ be given with $\parens*{\frac12 J''(w,x;\cdot)}\conjugate(h)<+\infty$. Then there exists a sequence $(h_k)$ with $h_k \in \frac12 \partial J''(w,x; \cdot)(z_k)$
for $z_k \in C_c(\Omega)$ such that $h_k \to h$ in $\MM(\Omega)$ and
\[
 \parens*{\frac12 J''(w,x;\cdot)}\conjugate(h_k) \to
 \parens*{\frac12 J''(w,x;\cdot)}\conjugate(h) .
\]
\end{lemma}
\begin{proof}
Due to \cref{theo312}, there exist densities $\hat v_1 \in L^1(\lambda)$, $\hat v_2 \in L^1(\HH^{d-1}|_\ZZ)$
	such that
	$h = \hat v_1 \lambda + \hat v_2 \HH^{d-1}|_{\ZZ}$
and
	\begin{equation}
		\label{eq:something_is_finite}
		\parens*{\frac12 J''(w;\cdot)}\conjugate(h)
		=
		\frac12  \int_{\ZZ} \frac{\abs{\nabla w}}{2a} \hat v_2^2 \d\HH^{d-1}
		+
		\int_\Omega \left(\frac12 j_0''(w;\cdot) \right)\conjugate(\hat v_1)  \d\lambda
		< \infty
		.
	\end{equation}
	From \cref{rem:conj},
	we get
	$\hat v_1 = 0$ $\lambda$-a.e.\ on $\set{j_0''(w; \sign(\hat v_1)) = 0}$.

	Next, we define the measure
	$\nu := \lambda + \HH^{d-1}|_{\ZZ}$.
	Note that $\nu$ is finite on compact subsets of $\Omega$,
	since $\ZZ$ is locally a $(d-1)$-dimensional submanifold.
	Thus, we can write $h = v \nu$, where $v \colon \Omega \to \R$
	is given by
	\begin{equation*}
		v(\omega) :=
		\begin{cases}
			\hat v_1(\omega) & \text{if } \omega \in \Omega \setminus \ZZ, \\
			\hat v_2(\omega) & \text{if } \omega \in \ZZ.
		\end{cases}
	\end{equation*}
	Similarly, we define a weight function $\rho \colon \Omega \to \R$
	by
	\begin{equation*}
		\rho(\omega) :=
		\begin{cases}
			\max j_0''(w(\omega); \pm 1) & \text{if } \omega \in \Omega \setminus \ZZ, \\
			\frac{2 a(\omega)}{\abs{\nabla w(\omega)}} & \text{if } \omega \in \ZZ.
		\end{cases}
	\end{equation*}
	Here,
	$\max j_0''(w(\omega); \pm 1)
	:=
	\max\set{j_0''(w(\omega); + 1), j_0''(w(\omega); - 1)}$.
	Finally,
	we set
	\begin{equation*}
		z(\omega) :=
		\begin{cases}
			0 & \text{if } \omega \in (\Omega \setminus \ZZ) \cap \set{ j_0''(w;\sign(v)) = 0}, \\
			\frac{1}{j_0''(w(\omega); \sign(v(\omega)))} v(\omega) & \text{if } \omega \in (\Omega \setminus \ZZ) \setminus \set{ j_0''(w;\sign(v)) = 0}, \\
			\frac{\abs{\nabla w(\omega)}}{2 a(\omega)} v(\omega) & \text{if } \omega \in \ZZ.
		\end{cases}
	\end{equation*}
	For $\omega \in \ZZ$
	we have
	\begin{equation*}
		\rho(\omega) \abs{z(\omega)}^p
		=
		\parens*{\frac{\abs{\nabla w(\omega)}}{2 a(\omega)}}^{p-1} \abs{\hat v_2(\omega)}^p
	\end{equation*}
	and this function belongs to
	$L^1(\HH^{d-1}|_\ZZ)$
	for all $p \in \set{1,2}$
	due to $\hat v_2 \in L^1(\HH^{d-1}|_{\ZZ})$
	and \eqref{eq:something_is_finite}.
	For $\omega \in (\Omega\setminus \ZZ) \setminus \set{j_0''(w;\sign(v))=0}$
	we utilize \cref{asm:structure_of_conjugate}
	to obtain
	$\max j_0''(w(\omega); \pm 1) \le C_j j_0''(w(\omega); \sign(v(\omega)))$.
	Thus,
	we find
	\begin{equation*}
		\rho(\omega) \abs{z(\omega)}^p
		=
		\frac{\max j_0''(w(\omega); \pm 1)}{j_0''(w(\omega); \sign(v(\omega)))^p}
		\abs{\hat v_1}^p
		\le
		C_j
		j_0''(w(\omega); \sign(v(\omega)))^{1-p}
		\abs{\hat v_1}^p
	\end{equation*}
	Due to $\hat v_1 \in L^1(\lambda)$
	and \eqref{eq:something_is_finite} in conjunction with \cref{rem:conj},
	these functions belong to $L^1(\lambda)$.
	Hence,
	$z \in L^1(\rho \nu) \cap L^2(\rho \nu)$.
	Moreover, $\rho \nu$ is finite on compact sets.
	Via \cref{lem:cts}
	we get a sequence $\seq{z_k} \subset C_c(\Omega)$
	such that
	$z_k \to z$
	in $L^p(\rho \nu)$ for all $p \in \set{1,2}$.
	Next, we set
	\begin{equation*}
		v_k(\omega) :=
		\begin{cases}
			\frac12 \partial j_0''(w(\omega); \cdot)(z_k(\omega))
			= j_0''(w(\omega); \sign(z_k(\omega))) z_k(\omega)
			& \text{if } \omega \in \Omega \setminus \ZZ, \\
			\frac{2 a(\omega)}{\abs{\nabla w(\omega)}} z_k(\omega) & \text{if } \omega \in \ZZ,
		\end{cases}
	\end{equation*}
	which leads to $h_k := v_k \nu \in \frac12 \partial J''(w,x;\cdot)(z_k)$,
	see \eqref{eq:subdif_j_pp}.
	From $z_k \to z$ in $L^1(\rho \nu)$
	we get
	\begin{align*}
		0
		&\leftarrow
		\norm{z_k - z}_{L^1(\rho \nu)}
		=
		\int_\Omega
		\abs{z_k - z}
		\max j_0''(w; \pm 1)
		\d\lambda
		+
		\int_\ZZ
		\abs{z_k - z}
		\frac{2 a}{\abs{\nabla w}}
		\d\HH^{d-1}
		\\
		&\ge
		\int_\Omega
		\abs{j_0''(w; \sign z_k) z_k - j_0''(w; \sign z) z}
		\d\lambda
		+
		\int_\ZZ
		\abs{z_k - z}
		\frac{2 a}{\abs{\nabla w}}
		\d\HH^{d-1}
		\\
		&=
		\int_\Omega
		\abs{v_k - v}
		\d\lambda
		+
		\int_\ZZ
		\abs{v_k - v}
		\d\HH^{d-1}
		\\
		&=
		\norm{h_k - h}_{\MM(\Omega)}
		.
	\end{align*}
	Similarly,
	from $z_k \to z$ in $L^2(\rho \nu)$
	we get
	\begin{align*}
		2\parens*{\frac12 J''(w,x;\cdot)}\conjugate(h_k)
		&=
		\int_{\set{j_0''(w;\sign(v_k)) \ne 0}}
		\frac{1}{j_0''(w;\sign(v_k))}
		\abs{v_k}^2
		\d\lambda
		+
		\int_\ZZ
		\frac{\abs{\nabla w}}{2 a}
		\abs{v_k}^2
		\d\HH^{d-1}
		\\&=
		\int_\Omega
		j_0''(w; \sign(z_k))\abs{z_k}^2
		\d\lambda
		+
		\int_\ZZ
		\frac{2 a}{\abs{\nabla w}}
		\abs{z_k}^2
		\d\HH^{d-1}
		\\&\to
		\int_\Omega
		j_0''(w; \sign(z))\abs{z}^2
		\d\lambda
		+
		\int_\ZZ
		\frac{2 a}{\abs{\nabla w}}
		\abs{z}^2
		\d\HH^{d-1}
		\\&=
		\int_\Omega
		\frac{1}{j_0''(w;\sign(v))}
		\abs{v}^2
		\d\lambda
		+
		\int_\ZZ
		\frac{\abs{\nabla w}}{2 a}
		\abs{v}^2
		\d\HH^{d-1}
		\\&=
		2\parens*{\frac12 J''(w,x;\cdot)}\conjugate(h)
		.
	\end{align*}
	For the convergence of the first integral,
	we used the dominated convergence theorem
	with the dominating function
	$\max j_0''(w;\pm1) g^2 \in L^1(\lambda)$,
	where $g$ comes from \cref{lem:cts}.
\end{proof}

\begin{theorem}\label{thm_G_striclty_twice}
Let \cref{asm:structure_of_conjugate} be satisfied.
	Let $x \in \dom(G)$, $w \in \partial G(x) \cap C_0(\Omega)$ be fixed
	such that
	$w$ belongs to
	$C^1(\Omega)$
	and satisfies
	\eqref{eq_ass_w_pw} and \eqref{eq_int_nablaw_z_finite}.
	Then $G$ is strictly twice epi-differentiable at $x$ for $w$ and
	\begin{equation*}
		\frac12 G''(x, w; h)
		=
		\parens*{\frac12 J''(w,x;\cdot)}\conjugate(h)
		\qquad
		\forall h\in \MM(\Omega).
	\end{equation*}
\end{theorem}
\begin{proof}
	We will apply  \cite[Lemma~3.2]{ChristofWachsmuth2019} with
	the setting
	$Q = \parens*{\frac12 J''(w,x;\cdot)}\conjugate$
	and
	$Z = C_c(\Omega)$.
	Condition (i) holds due to \cref{lem:lower_bound}.
	\cref{lem:upper_bound} yields condition (ii),
	while \cref{thm:recovery_seq} verifies condition (iii).
	Then \cite[Lemma~3.2]{ChristofWachsmuth2019} yields the claim.
\end{proof}

\subsection{Localized ascent lemma}
\label{subsec:ascent}
By dualizing the descent lemma (\cref{lem:descent}),
we obtain an ascent lemma.
\begin{lemma}
	\label{lem:ascent}
	We assume that $J$ is (Gâteaux) differentiable at $w$
	and that
	\eqref{eq:descent}
	is satisfied
	with $x = J'(w)$
	and $\Lambda, \eta > 0$.
	Then,
	\begin{equation}
		\label{eq:ascent}
		G(y)
		\ge
		G(x) + \dual{y - x}{w}
		+
		\frac{1}{2 \Lambda} \norm{y - x}_{L^1(\lambda)}^2
		\qquad
		\forall y \in L^1(\lambda), \norm{y - x}_{L^1(\lambda)} \le \eta \Lambda
		.
	\end{equation}
\end{lemma}
\begin{proof}
	Let $v \in C_0(\Omega)$ with $\norm{v - w}_{C_0(\Omega)} \le \eta$
	be arbitrary.
	We combine \eqref{eq:descent}
	with the Fenchel--Young inequality
	$\dual{y}{v} \le J(v) + G(y)$
	and
	$\dual{x}{w} = J(w) + G(x)$
	to obtain
	\begin{equation*}
		G(x) + \dual{y-x}{w} - G(y)
		\le
		\dual{x - y}{v-w} + \frac \Lambda 2 \norm{ v - w }_{C_0(\Omega)}^2
		.
	\end{equation*}
	Now, we can minimize the right-hand side w.r.t.\ $v$.
	Owing to the Hahn--Banach theorem,
	we have
	\begin{equation*}
		\frac{\norm{x - y}_{L^1(\lambda)}^2}{\Lambda}
		=
		\frac{\norm{x - y}_{\MM(\Omega)}^2}{\Lambda}
		=
		\sup\set*{
			\dual{y-x}{d}
			\given
			d \in C_0(\Omega),
			\norm{d}_{C_0(\Omega)}
			=
			\frac{\norm{x-y}_{\MM(\Omega)}}{\Lambda}
		}
		.
	\end{equation*}
	Thus,
	we can choose a sequence $\seq{v_k}$ with
	\begin{equation*}
		\norm{v_k - w}_{C_0(\Omega)} = \frac{\norm{x - y}_{L^1(\lambda)}}{\Lambda} \le \eta
		\quad\text{and}\quad
		\dual{x-y}{v_k-w} \to -\frac{\norm{x-y}_{\MM(\Omega)}^2}{\Lambda}.
	\end{equation*}
	This yields
	the desired
	\begin{equation*}
		G(x) + \dual{y-x}{w} - G(y)
		\le
		-\frac1{2\Lambda} \norm{ y - x }_{L^1(\lambda)}^2
		.
		\qedhere
	\end{equation*}
\end{proof}
\begin{remark}
	\label{rem:global_ascent}
	If $j_0$ is globally Lipschitz,
	one can show that \eqref{eq:ascent}
	holds for $\eta = \infty$,
	see \cref{rem:global_lipschitz}.
\end{remark}

\subsection{No-gap second-order optimality conditions}
Now we are ready
to
apply the theory of \cref{sec:second-order_conditions}
with the setting $X = \MM(\Omega)$ and $Y = C_0(\Omega)$.

We collect all assumptions on the structure of the optimization problem
in one place.

\begin{assumption}\label{asm_G}
	\begin{enumerate}
		\item The function $G \colon \MM(\Omega) \to \bar \R$ is defined by
	\[
	G(x):=
	\begin{cases}
	\int_\Omega g\parens*{ \frac{\d x}{\d\lambda}(\omega) } \d\lambda(\omega)
	&\text{if } x \ll \lambda, \\
	+\infty & \text{otherwise,}
	\end{cases}
	\qquad\forall x \in \MM(\Omega)
	\]
	where $g \colon \R \to \bar\R$ is convex and lower semicontinuous.

		\item
			The mapping $F \colon \dom(G) \to \R$
			satisfies
			\cref{asm:standing_assumption}~\ref{asm:standing_assumption:3}
			with $X=\MM(\Omega)$, $Y=C_0(\Omega)$, $\bar x\in \dom(G) \subset L^1(\lambda)$,  $\bar w=-F'(\bar x) \in C^1(\Omega) \cap C_0(\Omega)$.
			Further, we assume that $h \mapsto F''(\bar x) h^2$ is sequentially weak-$\star$ continuous.

	\item The conjugate of $g$ has the form
	\[
		g\conjugate (w)
	=
	j_0(w)
	+
	\sum_{i = 1}^m a_i \abs{ w - b_i }
	\qquad\forall w \in \R
	\]
    with $m \in \N$, $a_i > 0$, $i=1\dots m$,  and distinct numbers $b_i\in \R$, $i=1\dots m$.
    In addition, $j_0$ is in $C^1(\R)$ with locally Lipschitz continuous  $j_0'$  such that
 \[
  j_0''(w;z):=\lim_{t\searrow0} \frac2{t^2} (j_0(w+t z)- j_0(w) - t j_0'(w) z) \in \R
 \]
 exists for all $w,z\in \R$.
\item
	For each bounded set $W\subset \R$ there exists a constant $C_j \ge 1$ such that the smooth part $j_0$
	satisfies
	\begin{equation*}
		j_0''(w; -1) = 0
		\text{ or }
		j_0''(w; +1) = 0
		\text{ or }
		C_j^{-1} j_0''(w; -1)
		\le
		j_0''(w; +1)
		\le
		C_j j_0''(w; -1)
	\end{equation*}
	for all $w \in W$.

	\end{enumerate}
\end{assumption}

\begin{assumption}\label{asm_w}
	\begin{enumerate}
		\item The function $\bar w\in C^1(\Omega)$ satisfies
\begin{equation*}
 \{\bar w=b_i \}  \cap \{\nabla \bar w=0\}=\emptyset \qquad \forall i=1\dots m.
\end{equation*}
\item
	There exist $C, \eta > 0$ with
	\begin{equation*}
		\sum_{i = 1}^m \lambda\parens{\set{ \abs{ \bar w - b_i } \le \varepsilon }}
		\le
		C \varepsilon
		\qquad
		\forall \varepsilon \in [0,\eta].
	\end{equation*}
	\end{enumerate}

\end{assumption}

\begin{theorem}[strictly twice epi-differentiability of \texorpdfstring{$G$}{G}]
\label{thm_main_1}
Let \cref{asm_G,asm_w} be satisfied
and assume $\bar w \in \partial G(\bar x)$.
	Then $G$ is strictly twice epi-differentiable at $\bar x$ for $\bar w$.

Let $\mu \in \MM(\Omega)$.
If there exist densities $v_1 \in L^1(\lambda)$, $v_2 \in L^1(\HH^{d-1}|_\ZZ)$
	such that
	$\mu = v_1 \lambda + v_2 \HH^{d-1}|_{\ZZ}$
then
	\begin{equation*}
		\frac12 G''(\bar x, \bar w; \mu)
		=
		\frac12  \int_{\ZZ} \frac{\abs{\nabla \bar w}}{2a} v_2^2 \d\HH^{d-1}
		+
		\int_\Omega \left(\frac12 j_0''(\bar w;\cdot) \right)\conjugate(v_1)  \d\lambda \in [0,+\infty].
	\end{equation*}
Otherwise  $\frac12 G''(\bar x, \bar w; \mu)=+\infty$.

Here, $\ZZ:= \bigcup_{i=1}^m \set{\bar w=b_i}$, and $a:\ZZ\to \R$ is defined by $a(\omega) = a_i$ if $\bar w(\omega) = b_i$.
\end{theorem}
\begin{proof}
	The claim of strict twice epi-differentiability is \cref{thm_G_striclty_twice}
	as well as the identity
	$
	\frac12 G''(\bar x, \bar w; h)
		=
		\parens*{\frac12 J''(\bar w,\bar x;\cdot)}\conjugate(h)
		$.
		The form of the latter is given by \cref{theo312}.
\end{proof}

\begin{theorem}[no-gap SOC for integral functionals]
	\label{thm:SSC_wo2}
Let \cref{asm_G,asm_w} be satisfied.
Then
	\begin{equation}
		\label{eq:SOC_no_gap_integral_lang}
		F''(\bar x) h^2
		+
		G''(\bar x, \bar w; h)
		>
		0
		\qquad
		\forall h \in \MM(\Omega) \setminus \set{0}
	\end{equation}
if and only if
there are $c>0$, $\epsilon>0$ such that
\[
		F(x)+G(x)
		\ge
		F(\bar x) + G(\bar x) +  \frac{c}{2} \norm{x - \bar x}^2_{\MM(\Omega)}
		\qquad\forall x \in B_\varepsilon^{\MM(\Omega)}(\bar x).
\]
\end{theorem}
\begin{proof}
	The reverse implication follows directly from \cref{thm:SNC}.

	For the direct implication we are going to employ \cref{thm:SSC_wo}.
	By \cref{lem:Gpp_convex},
	we have
	$G''(\bar x, \bar w; h) = -\infty$
	for some $h \in X \setminus \set{0}$
	if $\bar w \not\in \partial G(\bar x)$.
	Thus, \eqref{eq:SOC_no_gap_integral_lang}
	implies
	$\bar w = -F'(\bar x) \in \partial G(\bar x)$.
	Now, we combine this with \cref{lem:descent,lem:ascent}
	to obtain the growth condition \eqref{eq:suf_nds}.
	Consequently,
	\cref{lem:suf_nds}
	yields that
	condition \eqref{eq:NDS} is fulfilled.
\end{proof}
Note that,
due to $G \equiv +\infty$ on $\MM(\Omega) \setminus L^1(\lambda)$,
the quadratic growth condition is equivalent to
\[
		F(x)+G(x)
		\ge
		F(\bar x) + G(\bar x) +  \frac{c}{2} \norm{x - \bar x}^2_{L^1(\lambda)}
		\qquad\forall x \in B_\varepsilon^{L^1(\lambda)}(\bar x).
\]

\begin{remark}
	\label{rem:subdif_not_nec}
	We mention that the knowledge of the precise structure of the subderivative $G''(\bar x, \bar w; \cdot)$
	is not necessary for the formulation and verification
	of \cref{thm:SSC_wo2}.
	Indeed, we only have to verify
	\eqref{eq:NDS} in the proof.
\end{remark}

\section{Applications}
\label{sec:app}
We are going to apply the abstract theory
to two non-smooth optimal control problems.
Both problems share the same smooth functional
\begin{equation*}
	F(u)
	:=
	\int_\Omega
	L(\cdot, y_u)
	\d\lambda
	\qquad\forall u \in L^1(\lambda),
\end{equation*}
where $y_u$ is the weak solution of the nonlinear PDE
\begin{align*}
	-\Delta y_u + a(\cdot, y_u) &= u \quad\text{in } \Omega, &
	y_u &= 0 \quad\text{on }\partial\Omega.
\end{align*}
Here, $\Omega \subset \R^d$ with $d \in \set{1,2,3}$ is open and bounded.
We require that the nonlinearities
$a,L \colon \Omega \times \R \to \R$
and $\Omega$
satisfy the usual assumptions,
see, e.g., \cite[Section~2.1]{CasasWachsmuthsBangBang},
such that
$F \colon L^\infty(\lambda) \to \R$
is twice continuously Fréchet differentiable with
\begin{align}
	F'(u)v &= \int_\Omega \varphi_u v \d\lambda \label{E2.5},
	\\
	F''(u)(v_1,v_2)
	&=
	\int_\Omega\bracks[\Big]{\frac{\partial^2L}{\partial y^2}(\cdot,y_u) - \frac{\partial^2a}{\partial y^2}(\cdot,y_u) \varphi_u} z_{u, v_1}z_{u, v_2}\d\lambda,
\label{E2.6}
\end{align}
where
the linearized states
$z_{u, v_i}$
and the
adjoint state $\varphi_u$
solve
the linear equations
\begin{align*}
	-\Delta z_{u, v_i} + \frac{\partial a}{\partial y}(\cdot, y_u) z_{u, v_i}
	&=
	\mrep{v_i}
	{\frac{\partial L}{\partial y}(\cdot, y_u)}
	\quad\text{in }\Omega,
	&
	z_{u, v_i} &= 0
	\quad\text{on } \partial\Omega,
	\\
	-\Delta
	\mrep{\varphi_u}{z_{u,v_i}}
	+ \frac{\partial a}{\partial y}(\cdot, y_u)
	\mrep{\varphi_u}{z_{u,v_i}}
	&=
	\frac{\partial L}{\partial y}(\cdot, y_u)
	\quad\text{in }\Omega,
	&
	\varphi_u &= 0
	\quad\text{on } \partial\Omega.
\end{align*}
In particular,
the adjoint state $\varphi_u$
is a representative for $F'(u)$.

As in \cref{sec:sod_integral_functionals_measures},
we use the spaces
$Y = C_0(\Omega)$,
$X = Y\dualspace = \MM(\Omega)$.
As usual, we will continue to denote the controls by $u$,
i.e., the symbols $x$ from \cref{sec:second-order_conditions}
will change to $u$.
In order to comply with \cref{asm:standing_assumption},
we have to extend the bilinear form $F''(u)$
(which is defined on $L^\infty(\lambda)^2$)
to $\MM(\Omega)^2$
and this is possible
by regularity theory for elliptic equations
due to $d \le 3$,
see, e.g., \cite[Section~2.5]{CasasWachsmuthsBangBang}.
This extended bilinear form is even sequentially weak-$\star$
continuous on $\MM(\Omega)^2$.
It is now easy to see that the differentiability of $F$
implies
that \cref{asm:standing_assumption}
is satisfied at all $\bar u \in L^\infty(\lambda)$
as long as $\dom(G)$ is bounded in $L^\infty(\lambda)$.
Thus,
\cref{asm:standing_assumption}
is satisfied.

\subsection{Application to a bang-off-bang problem}
\label{sec:app_bang}
First,
we consider an example with a bang-off-bang solution structure,
see \cite[Section~3]{Casas2012:1}.
This is accomplished by using the functional
\begin{equation*}
	G(u) = \int_\Omega g(u) \d\lambda
	\qquad\text{with}\qquad
	g(u) := \alpha \abs{u} + \delta_{[u_a,u_b]}(u)
	.
\end{equation*}
We assume that the constants satisfy $\alpha > 0$ and $u_a < 0 < u_b$.
The convex conjugate of $g$ is given by
\begin{align*}
	j(w) := g\conjugate(w)
	&=
	\max\set{ u_b (w - \alpha), u_a (w + \alpha), 0}
	\\&
	=
	\frac{\abs{u_a}}{2} \abs{w + \alpha}
	+
	\frac{u_b}{2} \abs{ w - \alpha}
	+
	\frac{u_b + u_a}{2}
	+
	\frac{u_a - u_b}{2}\alpha
	.
\end{align*}
It is easy to check that
\cref{asm_G}
is satisfied.
Next,
we choose a point $\bar u \in \dom(G)$
such that
$\bar w := -\bar\varphi := -\varphi_{\bar u}$
satisfies
\cref{asm_w}
with
$b_1 = -\alpha$, $b_2 = \alpha$.
We set $\ZZ := \set{\abs{\bar\varphi} = \alpha}$
and we define $u_{a,b} := \abs{u_a} \chi_{\set{\bar\varphi = -\alpha}} + u_b \chi_{\set{\bar\varphi = \alpha}}$.
Consequently,
\cref{thm:SSC_wo2}
yields that
\begin{subequations}
	\label{eq:SOC_bangbangbang}
	\begin{align}
		F'(\bar x) + \partial G(\bar x) &\ni 0
		\qquad\text{and}
		\\
		\label{eq:SOC_bangbangbang_2}
		F''(\bar u) \parens[\Big]{v \HH^{d-1}|_{\ZZ}}^2
		+
		\int_{\ZZ} \frac{\abs{\nabla\bar\varphi}}{u_{a,b}} v^2 \d\HH^{d-1}
		&>
		0
		\qquad\forall
		v \in L^1(\HH^{d-1}|_{\ZZ}) \setminus \set{0}
	\end{align}
\end{subequations}
is equivalent to the existence of $c,\epsilon > 0$
with
\[
		F(u)+G(u)
		\ge
		F(\bar u) + G(\bar u) +  \frac{c}{2} \norm{u - \bar u}^2_{L^1(\lambda)}
		\qquad\forall u \in B_\varepsilon^{L^1(\lambda)}(\bar u).
\]

We briefly compare our result with
\cite[Theorem~3.6]{Casas2012:1}.
Therein, the author does not require \cref{asm_w} on the structure of $\bar w$,
but he only obtains a quadratic growth w.r.t.\ the $L^2$-norms of the linearized state $z_{\bar u, u - \bar u}$
and the difference of the states $y_u-y_{\bar u}$.
Moreover,
our second-order condition
\eqref{eq:SOC_bangbangbang_2}
is different
since the second term is not present in
\cite[Theorem~3.6]{Casas2012:1}
and we do not need an enlarged critical cone
but we can work directly with measures on $\ZZ$,
similar to \cite[Section~2.5]{CasasWachsmuthsBangBang}.

The bang-bang case ($\alpha = 0$)
works similarly
and
we can reproduce
the results from
\cite[Example~6.14]{ChristofWachsmuth2017:1}.
In a similar spirit, control problems with the convexified multi-bang functional from \cite{ClasonKunisch2014}
can be handled.

\subsection{Application to a relaxed \texorpdfstring{$L^0$}{L0}-problem}
\label{sec:app_l0}

As a second application, we consider a problem with the following functional
\[
	\tilde G(u) = \int_\Omega \tilde g(u) \d\lambda
	\qquad\text{with}\qquad
	\tilde g(u) := \frac\alpha2 u^2 +\beta \abs{u}_0 + \delta_{[-\gamma,\gamma]}(u),
\]
with positive constants $\alpha, \beta, \gamma$,
and
\[
\abs{u}_0 :=
	\begin{cases}
		1 & \text{ if } u\ne 0, \\
		0 & \text{ if } u=0.
	\end{cases}
\]
Clearly, $\tilde g$ is not convex. As in
\cite{CasasWachsmuth2020}, we will study second-order conditions for problems with the convex envelope of $\tilde g$.
The convex envelope of $\tilde g$ is given by
\begin{equation*}
    g(u)
    =
    \begin{cases}
        \sqrt{2 \alpha \beta} \abs{u} & \text{if } \abs{u} \le \sqrt{2 \beta / \alpha}, \\
        \frac\alpha2 u^2 + \beta & \text{if } \abs{u} \in [\sqrt{2 \beta / \alpha}, \gamma], \\
        +\infty & \text{if } \abs{u} > \gamma.
    \end{cases}
\end{equation*}
In case $\gamma \le \sqrt{2 \beta / \alpha}$, $g$ coincides with the functional considered in the previous section.
Hence, we will focus on the case $\gamma > \sqrt{2 \beta / \alpha}$.

\begin{lemma}
	\label{lem:convex_conjugate_of_g}
	The convex conjugate of $g$ is given by
	\begin{equation*}
		j(w)
		=
		g\conjugate(w)
		=
		\begin{cases}
			0 & \text{if } \abs{w} \le \sqrt{2 \alpha \beta}, \\
			\frac{w^2}{2 \alpha} - \beta & \text{if } \abs{w} \in [\sqrt{2 \alpha \beta}, \alpha \gamma], \\
			\gamma \abs{w} - \frac\alpha2 \gamma^2 - \beta & \text{if } \abs{w} \ge \alpha \gamma.
		\end{cases}
	\end{equation*}
\end{lemma}
\begin{proof}
	We recall that
	the convex conjugate is defined via
	\begin{equation*}
		g\conjugate(w)
		=
		\sup\set*{ u w - g(u) \given u \in \R}.
	\end{equation*}
	Since the effective domain of $g$ is compact, the supremum is always attained.
	In case $\abs{w} \le \sqrt{2 \alpha \beta}$,
	the supremum is attained at $u = 0$ and, thus,
	$g\conjugate(w) = 0$.
	We consider the case $\abs{w} \in (\sqrt{2 \alpha \beta}, \alpha \gamma)$.
	For $u = w / \alpha$ we have
	$\abs{u} \in (\sqrt{2 \beta/\alpha}, \gamma)$
	and, thus,
	$0 = w - \alpha u = w - g'(u)$.
	This yields
	\begin{equation*}
		g\conjugate(w)
		=
		u w - g(u)
		=
		\frac{w^2}{\alpha} - \frac{\alpha}2 \frac{w^2}{\alpha^2} - \beta
		=
		\frac{w^2}{2 \alpha} - \beta
		.
	\end{equation*}
	Finally,
	for $\abs{w} \ge \alpha \gamma$,
	the supremum is attained at $u = \sign(w) \gamma$
	and the announced formula follows.
\end{proof}

\cref{asm_G} is satisfied with $m=2$, $b_1,b_2 = \pm \sqrt{2 \alpha \beta}$, $a_1 = a_2 = \frac12\sqrt{2 \beta/\alpha} $.
The corresponding function $j_0$ satisfies
\[
	j_0''(w;z) =
	\begin{cases}
		\frac1\alpha & \text{ if } \abs{w} = \sqrt{2 \alpha \beta} \text{ and } \sign (w) = \sign(z), \\
		\frac1\alpha & \text{ if } \abs{w} \in (\sqrt{2 \alpha \beta}, \alpha \gamma),\\
		\frac1\alpha & \text{ if } \abs{w} = \alpha \gamma  \text{ and } \sign (z) = -\sign(w), \\
		0 & \text{ otherwise}.
	\end{cases}
\]
Its conjugate is given by
	\begin{equation*}
		\left(\frac12 j_0''(w;\cdot)\right)\conjugate (v) :=
		\begin{cases}
			\frac\alpha2 v^2
			& \text{if } \abs{w} \in ( \sqrt{2 \alpha \beta}, \alpha \gamma),
			\\
			\frac\alpha2 v^2 + \delta_{(-\infty,0]}
			& \text{if } w = +\alpha \gamma
			\text{ or } w= -\sqrt{2 \alpha \beta} ,
			\\
			\frac\alpha2 v^2 + \delta_{[0,+\infty)}
			& \text{if } w = -\alpha \gamma
			\text{ or } w= +\sqrt{2 \alpha \beta} ,
			\\
			\delta_{\set{0}}(y) & \text{else}.
		\end{cases}
	\end{equation*}
Note that $\left(\frac12 j_0''(\bar \varphi;\cdot) \right)\conjugate(v_1) <+\infty$
implies sign conditions on $v_1$.

Now let a control $\bar u \in \dom(G)$ be given
such that
$\bar w := -\bar\varphi := -\varphi_{\bar u}$
satisfies
$\bar w \in \partial G(\bar u)$
and
\cref{asm_w}.
We set $\ZZ := \set{\abs{\bar\varphi} = \sqrt{2 \alpha \beta}}$.
Let $\mu \in \MM(\Omega)$ be given.
If there exist densities $v_1 \in L^1(\lambda)$, $v_2 \in L^1(\HH^{d-1}|_\ZZ)$
	such that
	$\mu = v_1 \lambda + v_2 \HH^{d-1}|_{\ZZ}$
	and
	\begin{equation}\label{eq_l0_v1_sign}
		v_1 \le 0  \quad\text{$\lambda$-a.e.\ on } \set{\bar\varphi = +\alpha \gamma},
		\qquad
		v_1 \ge 0  \quad\text{$\lambda$-a.e.\ on } \set{\bar\varphi = -\alpha \gamma},
	\end{equation}
then
	\begin{equation*}
		\frac12 G''(\bar x, \bar \varphi; \mu)
		=
		\frac12  \int_{\ZZ} \sqrt{\frac{\alpha}{2 \beta}} \abs{\nabla \bar \varphi} v_2^2 \d\HH^{d-1}
		+
		\frac\alpha2
		\int_{ \set{\abs{\bar\varphi} \in (\sqrt{2 \alpha \beta}, \alpha \gamma]} } v_1^2  \d\lambda \in [0,+\infty].
	\end{equation*}
Otherwise  $\frac12 G''(\bar x, \bar \varphi; \mu)=+\infty$.
Note that \eqref{eq_l0_v1_sign}
does not contain sign conditions on $\set{\abs{\bar\varphi} = \sqrt{2\alpha\beta}} = \ZZ$,
since this is a null set due to \cref{asm_w}.
Due to the same reason,
we do not integrate over this set in the second integral above.

Using this second-order derivative of $G$, we can formulate second-order necessary and sufficient optimality conditions.
We compare this to the results of \cite{CasasWachsmuth2020}. There,
the following second-order necessary condition \cite[Theorem 4.11]{CasasWachsmuth2020} was proven
\[
F''(\bar u) \parens{v_1 \lambda}^2
		+ G''(\bar x, \bar \varphi; v_1\lambda) \ge 0
\]
for all $v_1$ satisfying the sign conditions \eqref{eq_l0_v1_sign}.
In addition, second-order sufficient conditions were obtained with a second-order expression $\tilde G''$
which is strictly smaller than the second derivative obtained in this work.
Following the argument of \cite[Section 4.4]{CasasWachsmuth2020}, our no-gap conditions can be translated into no-gap conditions for problems
with the non-smooth functional $\tilde g$.

\section{Conclusions and outlook}
\label{sec:concl}
We have derived no-gap second order conditions
for optimal control problems
of the form \eqref{eq:problem},
in which $F$ is smooth and depends only on the state variable
and
in which $G$ is determined
by a not uniformly convex integrand.
Thus, the classical theory (e.g., \cite{CasasTroeltzsch2012,Casas2015})
cannot be utilized.
We have seen that it is fruitful to discuss these problems
in the space $\MM(\Omega) = C_0(\Omega)\dualspace$.
The required weak-$\star$ second subderivatives
can be obtained
by dualizing results for the preconjugate function.

Our second order condition is equivalent to a linear growth in the space $L^1(\lambda)$.
This follows automatically since we have applied
the results of
\cref{sec:second-order_conditions}
with the space $\MM(\Omega)$.
We expect that a better growth estimate can be obtained
by using a different space $X$ which is tailored
to the precise form of the functional $G$.
This is subject to future work.

\bibliographystyle{jnsao}
\bibliography{references}

\end{document}